\documentclass[11pt]{article}
\usepackage[utf8]{inputenc}
\usepackage{amsbsy,amsfonts,amsmath,amssymb,amsthm
,bbm,enumerate,euscript,graphicx,color, subcaption}

\usepackage{hyperref} 
\hypersetup{
    colorlinks=true,
    linkcolor=green,
    citecolor=green,
    urlcolor = magenta,
}

\def\R{\mathbb{R}}
\def\N{\mathbb{N}}
\def\cN{\mathbf{N}}
\def\Z{\mathbb{Z}}
\def\E{\mathbb{E}}

\def\V{\mathrm{Var}}
\def\C{\mathrm{Cov}}
\def\P{\mathbb{P}}

\def\p{\mathbf{p}}
\def\I{\mathbf{1}}
\def\D{\mathrm{d}}

\newtheorem{theo}{Theorem}
\newtheorem{coro}{Corollary}

\newtheorem{prop}{Proposition}
\newtheorem{lemm}{Lemma}

\theoremstyle{definition}
\newtheorem{defn}{Definition}
\newtheorem*{proo}{Proof}

\numberwithin{equation}{section}
\newtheorem*{rema}{Remark}

\title{\bf {Fitting regular point patterns with a hyperuniform perturbed lattice}}
\author{Daniela Flimmel 
\thanks{Corresponding author: \texttt { daniela.flimmel@matfyz.cuni.cz }(Daniela Flimmel)\\address: Sokolovská 83, 186 75 Prague, Czech Republic}
\and
\footnotesize{Department of Probability and Mathematical Statistics, Charles University}}

\date{\today}

\begin{document}
\maketitle

\begin{abstract}
We introduce a flexible methodology for modelling regular spatial point patterns using hyperuniform perturbed lattices. We show that, under suitable mixing conditions on the displacement field, lattices perturbed by stationary random fields are hyperuniform in arbitrary dimension. In particular, Gaussian perturbations with \textcolor{black}{absolutely summable covariances yield class-I hyperuniform point processes}. We further derive an explicit formula for the $K$-function of Gaussian models, which enables efficient parameter estimation via the minimum contrast method. The proposed framework provides a computationally tractable alternative to classical Gibbs models for repulsive data. The methodology is illustrated on three-dimensional data describing grain centers in a polycrystalline nickel-titanium alloy.
\end{abstract}

\vspace{0.2cm}


\vspace{0.5cm}

\section{Introduction}
Hyperuniformity, often described as a state of \textit{global order and local disorder}, refers to a unique characteristic of certain point processes (or structures in general) whose fluctuations in the number of points are suppressed compared to the Poisson point process. The notion was introduced and systematically developed in the physics literature by \cite{TS03, T16, T18}. For a mathematical treatment, we refer to the surveys \cite{GL17, C21, LR25} and the references therein. This property appears in many scientific fields. For instance, the eigenvalues of Gaussian random matrices were shown to be hyperuniform, as well as the zeros of Gaussian random polynomials. In biology, avian photoreceptor patterns \cite{JLH14} exhibit hyperuniformity, as do receptors in the immune system \cite{MBMW15}. In particular, the concept of hyperuniformity is relevant in the study of physical systems, offering insights into the organization of particles at large
scales despite their seemingly disordered local appearance (e.g. the order of the early Universe \cite{P93}).

\textcolor{black}{The focus of this paper is on \textit{regular} point patterns, that is, point sets exhibiting a higher degree of spatial uniformity than completely random configurations. Such regularity is often accompanied by some form of local repulsion between nearby points, although the two notions should not be identified.} \textcolor{black}{However, local repulsion alone does not imply hyperuniformity. For example, hard-core models and random sequential adsorption processes exhibit short-range exclusion while their fluctuation behaviour differs substantially from that of hyperuniform systems (see \cite{SPY07, DF24}). Thus, although these models are more regular than a Poisson point process at short distances, their number variance remains of standard volume order rather than displaying the anomalous suppression of fluctuations characteristic of hyperuniform systems.} Formally, this connection can be understood through asymptotically negative correlations between the numbers of points in two neighboring rectangular windows, a phenomenon described in detail in \cite{JS25} and closely related to notions of rigidity and suppressed fluctuations in point processes (see \cite{GL17}). 

Regular patterns arise in many natural settings, including cellular structures and crystalline arrangements in various materials. \textcolor{black}{One important mechanism contributing to spatial regularity is short-range exclusion between nearby points, which is frequently modeled using Gibbs point processes with repulsive energy functions.} While such models typically provide a good fit to observed data, the fluctuations in the number of points per unit volume remain Poisson-like, as noted in \cite{DF24}. However, if there is evidence of hyperuniformity, either from the physical nature of the data or from a given statistical test (see recent results \cite{KLH24, HGB23, MBL24}), then important questions arise. What model should we choose, and how can we efficiently generate large samples from such a model for further statistical inference? While some determinantal point
processes, sine-$\beta$ processes, or Coulomb gases are known to be hyperuniform, generating large samples of these processes remains computationally challenging. Until recently, the only way to limit the computation times was to start with a fixed hyperuniform structure such as $\Z^d$ and displace all the points independently. \textcolor{black}{Unfortunately, this construction preserves a substantial amount of the long-range order inherited from the underlying lattice, leading to atoms in the spectral measure (Bragg peaks); see \cite{BH24, Y21, KKT20, PS14, H95} for further discussion of diffraction properties and the persistence or attenuation of Bragg peaks under random perturbations. Although computationally attractive, this feature makes the model too restrictive for the nitinol data analysed in Section \ref{S_4}.}

A promising new alternative lies in dependent displacements of a lattice (\cite{G04, C21, KKT20, PS14}). It was shown in \cite{DFHL24} that the finite second moment of the common law of perturbations ensures that the resulting process is hyperuniform in dimensions $d=1,2$ regardless of the covariance structure among the perturbations. However, in dimensions $d\geq 3$, this is far from true. Even almost surely bounded perturbations of the higher dimensional lattice can produce non-hyperuniform structures.

There are several theoretical and practical motivations for considering hyperuniform perturbed lattices as baseline models for regular spatial data. From a statistical perspective, reduced fluctuation order can improve numerical stability in integration-based estimation procedures. In particular, slower growth of the number variance enhances the efficiency of Monte Carlo–type approximations of statistics of the form $\int f \D \mu$. Moreover, recent results in \cite{LRY24} show that, under suitable integrability conditions on the pair correlation function, any hyperuniform point process in dimensions $d\leq 2$ lies within finite Wasserstein distance of an $L^2$-perturbed lattice. Related transport-based and variance asymptotic results for stationary random measures can be found in \cite{KLLY25}. \textcolor{black}{Since perturbed lattices naturally preserve the regular structure of the underlying lattice, they provide an attractive and interpretable modeling framework for experimental data exhibiting regularity.}

The contributions of this paper are threefold. From a \emph{mathematical perspective}, we derive general sufficient conditions for class-I hyperuniformity of perturbed lattices in arbitrary dimension, formulated in terms of the $\alpha$-mixing coefficients of the perturbation field. Recently, the same conclusion was shown in \cite{KLLY25} under the assumptions on $\beta$-mixing coefficients.
In particular, we show that Gaussian perturbations with sufficiently fast decaying correlations satisfy these conditions. From a \emph{statistical perspective}, we propose a flexible parametric modeling framework for regular spatial data based on lattices perturbed by stationary Gaussian random fields. A key feature of this model is that the associated $K$-function admits an explicit formula, regardless of whether hyperuniformity holds. From a \emph{computational perspective}, the explicit form of the $K$-function enables fast parameter estimation via the minimum contrast method and allows efficient simulation without the numerical burden typically associated with determinantal or Coulomb-type models. The theoretical results are complemented by an analysis of three-dimensional spatial data obtained from synchrotron X-ray diffraction measurements of a polycrystalline nickel-titanium alloy. The study demonstrates that the choice of correlation structure in the perturbation field plays a crucial and delicate role in practice.

\textcolor{black}{The remainder of the paper is organized as follows. In Section \ref{S_main}, we present our main theoretical results, establishing sufficient conditions for hyperuniformity of perturbed lattices in terms of mixing properties of the displacement field and improving these results to Gaussian perturbations. Section \ref{S_background} introduces the necessary background on point processes, hyperuniformity, perturbed lattices, and develops the probabilistic tools used in the proofs. In Section \ref{S_Gaussian}, we focus on Gaussian perturbation fields and derive an explicit expression for the K-function of the resulting perturbed lattice model, which forms the basis for our statistical inference in Section \ref{S_4}. Finally, Section \ref{S_Proofs} contains the proofs of the theoretical results and auxiliary lemmas.}

\section{Main results}\label{S_main}
Let $\mathcal{L}$ be a lattice in $\R^d$ (e.g. $\Z^d$ or its affine transformation). By a perturbed lattice, we understand a stationary \textcolor{black}{random counting measure of the form
$$\Xi=\sum_{i \in \mathcal{L}}\delta_{i+p_i+U},$$ }
where $\mathbf{p}=\{p_i, i \in \mathcal{L}\}$ is a stationary random field in $\R^d$ and $U$ an independent uniform random variable (see Section \ref{S_2.4} for detailed treatment).

Our main results are formulated in terms of the decay of dependence within the perturbation field $\mathbf{p}$. To measure this dependence, we employ the strong mixing coefficients introduced by Rosenblatt and developed in the random-field setting by \cite{D94}. For two sub-$\sigma$-fields $\mathcal{A}_1, \mathcal{A}_2 \subset \mathcal{A}$, we denote
\begin{equation}\label{Eq_mixingOfAlg}
\alpha(\mathcal{A}_1, \mathcal{A}_2) = \textrm{sup}\{|\P(A \cap B)-\P(A)\P(B)|;\, A \in \mathcal{A}_1, B \in \mathcal{A}_2\}.
\end{equation} Specifically, for $C \subset \mathcal{L}$, we denote by $\sigma(p_C)$ the $\sigma-$field generated by $p_C:= \{p_i, i \in C\}$. 
\textcolor{black}{For $k,\ell\in \mathbb N\cup\{\infty\}$ and $n>0$, define 
\begin{equation}\label{Eq_mixingOfField}
\alpha_{k,\ell}(n)
=
\sup\Bigl\{
\alpha\bigl(\sigma( p_A),\sigma( p_B)\bigr);\,
A,B\subset\mathcal L,\,
|A|\le k,\,
|B|\le \ell,\,
d_{\mathcal L}(A,B)\ge n
\Bigr\}, 
\end{equation}
where
$$d_{\mathcal{L}}(A,B)=\inf\{|i-j|;\, i \in A, j \in B\}, \qquad  A, B \subset \mathcal{L}.$$}

\begin{theo}\label{thm1}
    Let $\{p_i, i \in \mathcal{L}\}$ be a strictly stationary random field such that there exist $s > 0, q \geq d$ with $1/s+1/q=1$, $\E |p_0|^q<\infty$ and
    $$\sum_{i \in \mathcal{L}}\alpha^{1/s}_{1,1}(|i|)<\infty.$$ Then the corresponding perturbed lattice is hyperuniform. If, moreover, the law of $p_0$ is symmetric with $\E|p_0|^{2q}<\infty$, then $\Xi$ is class-I hyperuniform.
\end{theo}

Several important classes of perturbation fields satisfy the assumptions of Theorem \ref{thm1}. We say that the random field $\mathbf{p}$ is $m$-dependent, if $p_i, p_j$ are independent random variables whenever $|i-j|>m$.

\begin{coro}\label{thm2}
    Let $\{p_i, i \in \mathcal{L}\}$ be a strictly stationary random field with almost surely bounded perturbations. If 
$$\sum_{i \in \mathcal{L}} \alpha_{1,1}(|i|)<\infty,$$
then the corresponding perturbed lattice is class-I hyperuniform.
\end{coro}

\begin{coro}\label{coro1}
    For some $m \in \N$, let $\{p_i, i \in \mathcal{L}\}$ be a strictly stationary $m$-dependent random field such that $\E|p_0|^{d}<\infty$. Then the corresponding perturbed lattice is hyperuniform.
 \end{coro} 

These results support the conclusion that sufficiently strong long-range dependencies are necessary to break hyperuniformity of the stationarized lattice. Namely, the counterexample of almost surely bounded perturbations resulting in a hyperfluctuating point process constructed in \cite{DFHL24} must satisfy $\sum\alpha_{1,1}(|i|)=\infty.$

\textcolor{black}{
While Theorem \ref{thm1} is formulated in terms of mixing coefficients, Gaussian random fields admit a more explicit characterization through their covariance structure. The next result shows that, in the Gaussian setting, hyperuniformity follows from a simple summability condition on the covariance matrices.
\begin{theo}\label{Coro2}
Let $\{p_i, i \in \mathcal{L}\}$ be a strictly stationary Gaussian random field, such that $p_0\sim \mathcal{N}_d(\mu, \Sigma)$ for some $\mu \in \R^d$ and positive definite matrix $\Sigma \in \R^{d\times d}$. Assume, moreover, that
$$\sum_{i \in \mathcal{L}}\|\C(p_0,p_i)\|<\infty.$$
Then the corresponding perturbed lattice is class-I hyperuniform.
\end{theo}}
The proofs are postponed to Section \ref{S_Proofs} after we build up the necessary tools. 

\textcolor{black}{\section{Theoretical background}\label{S_background}}
\subsection{General notation}
We assume that all random variables are defined on a common probability space $(\Omega,\mathcal A,\mathbb P)$, which remains fixed throughout the paper. \textcolor{black}{For random vectors $X\in\mathbb R^q$ and $Y\in\mathbb R^p$, we write $\E X \in \R^q$ for the expectation of $X$ and $\C(X,Y)
:= \mathbb E\left[(X-\mathbb E X)(Y-\mathbb E Y)^\top\right]
\in\mathbb R^{q\times p}$ for the covariance matrix of $X$ and $Y$. In particular, $\V(X):=\C(X,X).
$ Scalars are identified with $1\times1$ matrices, so the same notation applies to real-valued random variables.} In addition, if two random vectors $X, Y$ are equal in distribution, we write $X\stackrel{D}{=}Y$. 

To describe spatial point processes, i.e. point processes on $\R^d$, $d\geq 1$, we denote by $\mathcal{B}^d$ the Borel $\sigma-$algebra on $\R^d$ and $\mathcal{B}_b^d$ the system of all bounded Borel sets. Moreover, $\lambda^d$ denotes the Lebesgue measure on $(\R^d, \mathcal{B}^d)$, $|\cdot|$ is the Euclidean norm, \textcolor{black}{$\|\cdot\|$ is the matrix norm} and $\#(A)$ stands for the cardinality of a set $A$. For any set $B \in \mathcal{B}^d$, we denote by $B^C:= \R^d \setminus B$ its complement. For the purpose of counting points in balls, we denote by $B(x,r)$ the closed Euclidean ball centered at $x$ with radius $r>0$. In particular, we use the abbreviation $B_r$ when $x=0$.

\subsection{Point processes and second-order measures}
At this point, we recall some terminology related to point processes. As usual, let $\cN$ be the space of \textcolor{black}{all locally finite integer-valued measures on $(\R^d, \mathcal{B}^d)$, that is,
$$
\cN:=
\left\lbrace
\mu=\sum_i n_i\delta_{x_i};
\, n_i\in\N,\ x_i\in\R^d,
\mu(B)<\infty,\ \forall B\in\mathcal{B}_b^d
\right\rbrace.
$$
The space $\cN$ is equipped with the smallest $\sigma$-algebra $\mathcal{N}$ making the mappings
$\mu\mapsto \mu(B),\, B\in\mathcal{B}^d$
measurable. The elements of $\cN$ are called \emph{counting measures}. We denote by $\cN_f\subset\cN$ the subset of finite counting measures.}

\textcolor{black}{A \textit{point process} is a random element $\mu$ defined on $(\Omega,\mathcal{A},\P)$ with values in $(\cN,\mathcal{N})$. A point process is called \textit{simple} if
$$
\P(\mu({x})\leq 1 \text{ for all } x\in\R^d)=1.
$$
In this case, $\mu$ can be identified with the random set of its atoms and we may write $x\in\mu$ whenever $\mu({x})=1$. Throughout the paper, perturbed lattices are understood as random counting measures. Simplicity is not required for the theoretical results, although the Gaussian perturbation models considered in Section \ref{S_Gaussian} are simple almost surely due to the absolute continuity of the Gaussian distribution.}

Finally, we call a point process $\mu$ \textit{stationary} if its distribution is invariant with respect to shifts, meaning that $\mu \stackrel{D}{=}\mu+x$ for all $x \in \R^d$.

Hyperuniformity of a point process is closely tied to its second-order measures. 
In the point process literature, these characteristics are typically described in terms of (factorial) moment and cumulant measures. By contrast, the hyperuniformity literature (\cite{C21, LR25}) most often formulates second-order information using the covariance measure of the process. In this section, we briefly review the relevant definitions and clarify the relationships among these objects. 

\textcolor{black}{
\begin{defn}[$m$-th moment and factorial moment measures]
Let $m\in\N$ and assume that $\E[\mu(B)^m]<\infty$
for all $B\in\mathcal B_b^d$.
The \textit{$m$-th moment measure} $\mu^{[m]}$ is the Borel measure on $(\R^d)^m$ defined by
$$\mu^{[m]}(B_1\times\cdots\times B_m)
:=
\E\mu^{\otimes m}(B_1\times\cdots\times B_m),$$
where $\mu^{\otimes m}$ denotes the $m$-fold product measure of the random counting measure $\mu$.
The \textit{$m$-th factorial moment measure} $\mu^{(m)}$ is defined as the expectation of the $m$-th factorial power measure of $\mu$ (see \cite[Chapter~8]{DVJ08}).
If $\mu$ is simple, then
\[
\mu^{(m)}(B_1\times\cdots\times B_m)
=
\E\sum_{x_1,\ldots,x_m\in\mu}^{\neq}
\mathbf 1(x_1\in B_1,\ldots,x_m\in B_m).
\]
 Here, the superscript $\neq$ indicates that the sum is taken over all ordered $m$-tuples of distinct points of $\mu$, i.e., diagonal terms are excluded.\end{defn}
}

\begin{defn}[$m$-th cumulant and factorial cumulant measures] Let $m \in \N$ and assume that $\E[\mu(B)^m]<\infty$ for all $B\in \mathcal{B}^d_b$. Then we define the $m$-th cumulant measure $\gamma^{[m]}$ and factorial cumulant measure $\gamma^{(m)}$ of a point process $\mu$ as measures on $(\R^d)^m$ given by
 $$ \gamma^{[m]}(B_1\times\cdots\times B_m)=\sum_{j=1}^m (-1)^{j-1}(j-1)!\sum\limits_{\substack{K_1\cup\cdots\cup K_j = \{1,\ldots,m\}\\ K_a\cap K_b=\emptyset,\, a\neq b}} \prod_{r=1}^j \mu^{[\#(K_r)]}(\times_{s \in K_r}B_s),$$ 
 $$ \gamma^{(m)}(B_1\times\cdots\times B_m)=\sum_{j=1}^m (-1)^{j-1}(j-1)!\sum\limits_{\substack{K_1\cup\cdots\cup K_j = \{1,\ldots,m\}\\ K_a\cap K_b=\emptyset,\, a\neq b}} \prod_{r=1}^j \mu^{(\#(K_r))}(\times_{s \in K_r}B_s).$$  
\end{defn}
Cumulant measures provide a convenient framework for quantifying and analyzing dependence between spatially separated regions of a point process. 
Their defining relation to factorial moment measures mirrors the classical relationship between mixed moments and mixed cumulants in probability theory.

For $m=1$, the moment, factorial moment, and cumulant measures all coincide and reduce to the measure $\mu^{[1]}$, called the \textit{intensity measure}. 
If $\mu$ is stationary, then $\mu^{[1]}$ is proportional to Lebesgue measure, i.e.,
$\mu^{[1]} = \rho \lambda^d$ for some $\rho \geq 0$, where $\rho$ is referred to as the \textit{intensity}.

For $m=2$, the cumulant measure $\gamma^{[2]}$ is called the \textit{covariance measure}, since it is the measure on $\R^d \times \R^d$ satisfying
\textcolor{black}{$$
\gamma^{[2]}(B_1\times B_2)
= \C\bigl(\mu(B_1), \mu (B_2)\bigr).
$$}
Assuming stationarity of $\mu$, second-order measures admit a disintegration that leads to their reduced counterparts. 
In particular, the \textit{reduced covariance measure} $\gamma^{[2]}_{\mathrm{red}}$ is the measure on $\R^d$ satisfying
$$
\gamma^{[2]}(\D x, \D y)
= \rho \gamma^{[2]}_{\mathrm{red}}(\D (y - x))\D x.
$$
The measure $\gamma^{[2]}_{\mathrm{red}}$ is a signed positive semidefinite measure in the sense that
$$
\int f(x)\,\overline{f(x+y)}\,\gamma^{[2]}_{\mathrm{red}}(\D y)\,\D x
= \V\!\left(\textcolor{black}{\int f(x) \mu (\D x)}\right) \geq 0
$$
for all bounded, compactly supported functions $f:\R^d \to \mathbb{C}$.
By Bochner’s theorem, there exists a non-negative measure $\Gamma$ on $\R^d$ such that
$$
\V\left(\textcolor{black}{\int f(x) \mu (\D x)}\right)
= \frac{1}{(2\pi)^d} \int |\hat{f}(k)|^2 \, \Gamma(\D k),
$$
where
$$
\hat{f}(k) :=\mathcal{F}[f](k)=  \int_{\R^d} f(x) e^{-i k \cdot x}  \D x
$$
denotes the Fourier transform of $f$.
The measure $\Gamma$ is called \textit{Bartlett’s spectral measure}, and, when it is absolutely continuous with respect to Lebesgue measure, its density $S$ is referred to as the \textit{structure factor}, i.e. $S(k):= \int e^{-i k \cdot x}\gamma^{[2]}_{red}(\D x)$ if this expression makes sense.

More commonly in the mathematical literature, the structure factor is expressed in terms of the reduced factorial cumulant measure $\gamma^{(2)}_{\mathrm{red}}$. 
This is based on the relation between the reduced cumulant and factorial cumulant measures
\begin{equation}\label{EQ_RedCumMeasures}
\gamma^{[2]}_{red}(\D x)=  \delta_0 (\D x) + \rho \gamma^{(2)}_{red}(\D x),
\end{equation}
where $\delta_0$ is the Dirac delta at $0$.
If $\gamma^{(2)}_{\mathrm{red}}$ has finite total variation, that is, if
$$|\gamma_{red}^{(2)}|(\R^d)<\infty,$$
then the structure factor admits the representation
\begin{equation}\label{EQ_StructureFactor}
S(k)
= 1 + \rho \int_{\R^d} e^{-i k \cdot x}\, \gamma^{(2)}_{\mathrm{red}}(\D x),
\quad k \in \R^d
\end{equation}
which describes the spectral density of fluctuations at frequency $k$.

Furthermore, disintegration of the second factorial moment measure $\mu^{(2)}$ leads to the reduced second factorial moment measure $\mu^{(2)}_{\mathrm{red}}$, also known as the \textit{pair correlation measure}. 
It is characterized by
\textcolor{black}{\begin{equation*}
\mu^{(2)}_{\mathrm{red}}(B)
= \E \int \I\{x \in [0,1]^d;\, y-x \in B\}\mu^{(2)}(\D x, \D y),
\end{equation*}}
see \cite[Chapter 8]{LP17}. Associated with $\mu^{(2)}_{\mathrm{red}}$ is \textit{Ripley’s $K$-function}, defined as the cumulative mass of $\mu^{(2)}_{\mathrm{red}}$ over a ball of radius $r \geq 0$,
\begin{equation}\label{D_K}
K(r) = \mu^{(2)}_{\mathrm{red}}(B_r).
\end{equation}
The $K$-function measures clustering or regularity at different scales. In practice, this is visualized by centered \textit{Besag's $L$-function} $L(r)-r$ given by
\begin{equation}\label{D_L}
    L(r):= \left(\frac{K(r)}{\lambda^d(B(0,1))}\right)^{1/d}.
\end{equation}
Then $L(r)-r>0$ indicates clustering at scale $r$ while $L(r)-r<0$ shows regularity.  
If, moreover, $\mu^{(2)}_{\mathrm{red}}$ is absolutely continuous with respect to $\lambda^d$, then its density $g$ is called the \textit{pair correlation function}. 
If $g-1\in L^1$, then \eqref{EQ_StructureFactor} reduces to the familiar expression (\cite{C21,MBL24})
$$
S(k)
= 1 + \rho \int_{\R^d} (g(x) - 1)\, e^{-i k \cdot x} \, \D x.
$$
The latter is based on the relation
\begin{equation*}
    \gamma_{red}^{(2)}(\D x) = \mu_{red}^{(2)}(\D x)- \D x.
\end{equation*}

\subsection{Generalized structure factor}\label{S_GeneralSF}
In the previous section, we recalled the definition of the structure factor for point processes such that $|\gamma_{red}^{(2)}|(\R^d)<\infty$. This assumption is violated for perturbed lattices (see Section \ref{S_2.4} below) and therefore, the structure factor in Definition \ref{EQ_StructureFactor} is not well defined as a function. \textcolor{black}{Following the classical spectral framework of \cite[Chapter 8]{DVJ08} and its recent use in the hyperuniformity literature \cite{C21, LR25}, we employ a distributional formulation of Bartlett's spectral measure. In particular, the singular mass at $0$ is analyzed through suitable smooth approximations, allowing us to characterize hyperuniformity in terms of the low-frequency behavior of the spectrum.}

Let $\mathbb{S}$ denote the space of all Schwartz functions (i.e. smooth functions such that all their derivatives, including the functions themselves, decay faster than every polynomial at infinity). The dual space $\mathbb{S}^*$ to $\mathbb{S}$ is the space of continuous linear mappings on $\mathbb{S}$ called \textit{distributions}. By convention, the action of a distribution $T \in \mathbb{S}^*$ on $\mathbb{S}$ is represented by a `duality pairing'
$$ \langle T,f\rangle:= T(f), \quad f \in \mathbb{S}.$$
Fourier transform of a distribution $T$ is again a distribution $\hat{T}$ given by
    $$\langle\hat{T},f\rangle =\langle T,\hat{f}\rangle,\quad f \in \mathbb{S}.$$
If $\nu$ is a measure on $\R^d$ such that $T_\nu(f)= \int f \D \nu$ defines a distribution. Then we write $\langle\nu,f\rangle=T_\nu(f)$ and
$$\langle\hat{\nu}, f\rangle = \int \hat{f}\D \nu.$$
Note that here, the choice $\nu = \gamma^{(2)}_{red}$ (if it exists) always induces a well-defined distribution since $\gamma^{(2)}_{red}$ is a linear combination of a positive semidefinite measure and Dirac delta, see \eqref{EQ_RedCumMeasures}.

For $r>0$, let
 $$j_r(x)=\frac{1} {\kappa_d (2\pi)^d}\frac{(J_{d/2}(r|x|)^2}{|x|^d},$$
where $J_{d/2}$ is the Bessel function of the first kind and $\kappa_d$ is a constant given in Eq. (58) of \cite{T18}. Then $j_r \in \mathbb{S}$, $\int j_r(x)\D x =1$ and
$$\hat{j_r} = \frac{\I_{B_r}\star \I_{B_r}}{\lambda^d(B_r)},$$
where $\star$ denotes the convolution of two functions. 
With regard to \eqref{EQ_StructureFactor}, we may define the structure factor as a distribution $S \in \mathbb{S}^*$ satisfying
\begin{equation}\label{EQ_generalSF}
     \langle S,j_r\rangle = 1 + \langle {\gamma}_{red}^{(2)},\hat{j_r}\rangle.
\end{equation}
\textcolor{black}{As a consequence of the Plancherel formula applied to the ball indicator, we obtain}
   $$\frac{\V (\textcolor{black}{\mu (B_r)})}{\lambda^d(B_r)}= \langle S, j_r \rangle.$$

\subsection{Hyperuniformity of a point process}

\begin{defn}[Hyperuniform point process]
A stationary point process $\mu$ is called \textit{hyperuniform}, if 
\begin{equation}\label{EQ_hyperuniformity}
    \sigma(r):= \frac{\V (\textcolor{black}{\mu(B_r)})}{\lambda^d(B_r)} \xrightarrow{r \to \infty} 0.
\end{equation}
\end{defn}
We can classify the hyperuniformity by the speed of the convergence of $\sigma$. In the seminal paper \cite{T18}, it was proposed to call a process $\mu$
\begin{itemize}
    \item \textit{class-I hyperuniform} if $\sigma(r) \sim c r^{-1}$
    \item \textit{class-II hyperuniform} if $\sigma(r) \sim c r^{-1} \log(r)$
    \item \textit{class-III hyperuniform} if $\sigma(r) \sim c r^{-\tau}$, $\tau \in (0,1)$.
\end{itemize}
The papers \cite{T18, C21} offer a wide range of examples from a mathematical and physical perspective. However, as mentioned in \cite{DFHL24}, the classification is not exhaustive and does not contain all hyperuniform point processes.

If the reduced covariance measure of $\mu$ has finite total variation, then having \eqref{EQ_StructureFactor}, hyperuniformity translates to $S(0)=0$. Moreover, we may assume that the structure factor follows
a power-law behavior near zero:
\begin{equation}\label{EQ_SFF}
    S(k) = O(|k|^{\tau}).
\end{equation}
If \eqref{EQ_SFF} holds, we say that the \textit{hyperuniformity exponent} is at most $\tau$. Hyperuniformity then occurs with $\tau >0$. There is a correspondence between the hyperuniformity exponent and \eqref{EQ_hyperuniformity}. Namely, $\tau \in (0,1)$ if $\V(\textcolor{black}{\mu(B_r)})= O(r^{d-\tau})$ and $\tau>1$ corresponds to $\V(\textcolor{black}{\mu(B_r)})= O(r^{d-1})$, i.e. class-I hyperuniformity. The structure factor exponent of an independently perturbed lattice, where the perturbations follow a symmetric law, cannot exceed $2$, see Corollary 2.1 of \cite{LR25}. For Gaussian perturbations, \cite{Y22} shows that the exponent is exactly $2$.

Generally, when $|\gamma_{red}^{(2)}|(\R^d)=\infty$ (which is the case of perturbed lattices), we conclude that the process $\mu$ is hyperuniform if
$$\langle S, j_r\rangle\xrightarrow{r \to \infty} 0 $$
in the sense of \eqref{EQ_generalSF}.

\subsection{Perturbed lattices}\label{S_2.4}
Our core concept is the deterministic lattice $\mathcal{L}$ in $\R^d$. The lattice $\mathcal{L}$ is defined as a subgroup of $\R^
d$ generated by linear combinations with integer coefficients of the vectors of some basis $(a_1,\ldots, a_d)$ of $\R^d$. An elementary example is
the lattice $\Z^d$. Any other lattice can be obtained as an affine transformation of $\Z^d$. The fundamental domain of $\mathcal{L}$ is denoted by $\mathcal{D}$ and defined by
$$\mathcal{D}:= \left\lbrace \sum_{i=1}^d t_i a_i; \,(t_1,\ldots,t_d) \in [0,1)^d\right\rbrace.$$ 
For $\mathcal{L}=\Z^d$, this is simply $\mathcal{D}=[0,1]^d$.

\begin{defn}[Perturbed lattice]\label{D_PerturbedLattice}
    Let $\p=\{p_i, i \in \mathcal{L}\}$ be a strictly stationary random field in $\R^d$, i.e.
    $$\{p_i, i \in \Lambda\} \stackrel{D}{=} \{p_{i+j}, i \in \Lambda\}$$
    for all $\Lambda \subset \mathcal{L}$ and $j \in \mathcal{L}$.
    Let $U$ be a uniform random variable $U$ on $\mathcal{D}$ independent of $\mathbf{p}$.
    By a \textit{perturbed lattice}, we understand a stationary point process $\Xi$ on $\R^d$ defined by
    \textcolor{black}{$$\Xi : =  \sum_{i \in \mathcal{L}}\delta_{i+p_i+U}.$$}
Specifically, if $\p$ consists of independent random variables, we speak about \textit{independently perturbed lattice}. If for all $i\in \mathcal{L}$, $p_i=0$ almost surely, then $\Xi$ is called \textit{stationarized lattice}.
\end{defn}

The second-order measures of a perturbed lattice $\Xi$ are well defined if $\E[\textcolor{black}{\Xi(B_r)^2}]<\infty$ for all $r>0$. By Lemma 2.6. of \cite{DFHL24}, this is true provided that $\E|p_0|^d<\infty.$ Under this condition, the reduced factorial moment measure $\mu^{(2)}$ is a locally finite measure given by
\begin{equation}
   \label{alpha2PSL}
  \mu_{red}^{(2)}(B) = \E \sum_{i \in \mathcal{L} \setminus \{0\}} \I{\{i+p_i-p_0 \in B\}}
   \end{equation}
for any Borel $B\subset \R^d$ (cf. Proposition 2.7 in \cite{DFHL24}). As a result, the structure factor of a perturbed lattice cannot be defined in the sense of \eqref{EQ_StructureFactor} since
$$\gamma^{(2)}_{red} (B)= \E \sum_{j \in \mathcal{L}\setminus \{0\}}\I\{j+p_j-p_0 \in B\}- \lambda^d(B), \quad B \in \mathcal{B}^d_b,$$
can have infinite total variation.

\textcolor{black}{We conclude this section by discussing the classical examples of stationarized and independently perturbed lattices. The results presented below are not new. The hyperuniformity of stationarized lattices was established in \cite{BH24}. In the case of independent Gaussian perturbations, the corresponding statements follow from \cite{ST06} and \cite{Y21}, and the arguments of the latter paper can be extended to perturbations with more general distributions. Independent perturbations of lattices have also been studied extensively in the physics literature; see, for example, \cite{G04,H95}. Finally, in dimensions $d\leq 2$, the conclusions below may also be deduced from the more general results of \cite{DFHL24}.} Nevertheless, these results are non-trivial, and a general, self-contained proof does not appear to be available in the literature. Since the proofs of our main results rely crucially on these facts, we include complete proofs in Section \ref{S_Proofs} for the sake of a clean presentation. 
\begin{lemm}[Hyperuniformity of the stationarized lattice]\label{L_statL}
Let $\Xi$ be a stationarized lattice, that is the point process
    $$\Xi : = \textcolor{black}{\sum_{i \in \mathcal{L}}\delta_{i+U}}.$$
Then for $r>0$
$$\langle S, j_r\rangle=  (2\pi)^d\sum_{x \in \mathcal{L}^*\setminus \{0\}} j_r\left(2\pi x\right).$$
and consequently, $\Xi$ is class-I hyperuniform.
\end{lemm}
\begin{lemm}[Hyperuniformity of independently perturbed lattices]\label{L_IPL}
    Assume that $\Xi$ is an independently perturbed lattice, that is the point process 
    $$\Xi:= \textcolor{black}{\sum_{i \in \mathcal{L}}\delta_{i+p_i+U}},$$
where $U, \{p_i, i \in \mathcal{L}\}$ are independent random variables. Denote by $\varphi_{p_0}$ the characteristic function of $p_0$. Then 
$$\langle S, j_r\rangle =  1- \int |\varphi_{p_0} (x)|^2 j_r(x) \D x +  (2\pi)^d \sum_{x \in \mathcal{L}^*\setminus\{0\}} |\varphi_{p_0}(2\pi x)|^2 j_r(2\pi x)$$
and $\Xi$ is hyperuniform. 
\end{lemm}
\textcolor{black}{If, moreover, $p_0$ has a finite first moment, then the independently perturbed lattice is class-I hyperuniform, see \cite{GS75}.}
If the law of $p_0$ is symmetric and not almost surely zero, the hyperuniformity exponent of an independently perturbed lattice is at most $2$ (cf. Corollary 2.1 of \cite{LR25}).

\section{$K-$function of the Gaussian perturbation model}\label{S_Gaussian}
\textcolor{black}{
The hyperuniformity results of Section \ref{S_main} apply to broad classes of perturbation fields. For statistical inference, however, one needs models with tractable second-order characteristics. We therefore focus on Gaussian perturbation fields, for which the $K$-function (and consequently the Besag's $L$-function) can be computed explicitly and used for efficient parameter estimation. Note that the results in this section apply beyond the hyperuniform models. For simplicity, we restrict ourselves to $\mathcal{L}=\Z^d$.}

\textcolor{black}{
 Let $\{p_i\}_{i\in\mathbb Z^d}$ be a strictly stationary Gaussian random field in $\mathbb R^d$ with mean $\mu\in\mathbb R^d$. We denote \[ \Sigma:=\V(p_0), \qquad C_i:=\C(p_i,p_0), \qquad i\in\mathbb Z^d. \] Throughout this section we consider the isotropic case, that is \[ \Sigma=\sigma^2 I_d, \qquad C_i=\rho_i \Sigma, \] for some $\sigma^2>0$ and correlation coefficients $\rho_i\in[-1,1]$. In particular, \[ \V(p_i-p_0) = 2(\Sigma-C_i) = 2(1-\rho_i)\Sigma = 2\sigma^2(1-\rho_i)I_d. \]}

\begin{prop}\label{L_2} Assume that \[ \Sigma=\sigma^2 I_d, \qquad C_i=\rho_i\Sigma, \] with $\sigma^2>0$, and let \[ \mathcal A := \{i\in\mathbb Z^d\setminus\{0\}; \,\rho_i=1\}. \] Then the $K$-function of the corresponding perturbed lattice $\Xi$ is 
\textcolor{black}{
\[ K(r) = \#(\mathcal A\cap B(0,r)) + \sum_{i\in\mathbb Z^d\setminus(\mathcal A\cup\{0\})} P_d(r^2/\kappa_i,|i|^2/\kappa_i), \qquad r\ge0, \]}
where \[ \kappa_i:=2\sigma^2(1-\rho_i), \] and $P_d(\cdot,\eta)$ denotes the cumulative distribution function of the non-central chi-squared distribution with $d$ degrees of freedom and non-centrality parameter $\eta$. \end{prop}

\begin{proo} Using \eqref{alpha2PSL}, we obtain \
\begin{align*} K(r) &= \mu^{(2)}_{\mathrm{red}}(B(0,r))\\ &= \E\sum_{i\in\mathbb Z^d\setminus\{0\}} \mathbf 1_{\{i+p_i-p_0\in B(0,r)\}}\\ &= \sum_{i\in\mathbb Z^d\setminus\{0\}} \P(i+p_i-p_0\in B(0,r))\\ &= \sum_{i\in\mathbb Z^d\setminus\{0\}} \P(|i+p_i-p_0|\le r), 
\end{align*} 
where the exchange of expectation and summation is justified by the non-negativity of the summands. For $i\in\mathbb Z^d\setminus\{0\}$, let us define $Y_i:=p_i-p_0.$ Then \[ \V(Y_i) = 2(\Sigma-C_i) = 2\sigma^2(1-\rho_i)I_d. \] 
\textcolor{black}{If $i\in\mathcal A$, then $\rho_i=1$ and $\V(Y_i)=0.$ Consequently, $Y_i=0$ almost surely and 
\begin{equation}\label{Eq_iInA}
 \P(|i+Y_i|\le r) = \mathbf 1_{\{|i|\le r\}}. 
\end{equation}}

\textcolor{black}{Assume now that $i\notin\mathcal A$ and set $\kappa_i:=2\sigma^2(1-\rho_i).$ Then $Y_i\sim \mathcal N_d(0,\kappa_i I_d).$ Further, if we define $Z_i:=\kappa_i^{-1/2}(i+Y_i)$, then $ Z_i \sim \mathcal N_d(\kappa_i^{-1/2}i,I_d).$ Consequently, $|Z_i|^2$ has a non-central chi-squared distribution with $d$ degrees of freedom and non-centrality parameter  $\eta_i = \kappa_i^{-1}|i|^2.$ Therefore 
\begin{align} 
\P(|i+Y_i|\le r) &= \P(|Z_i|^2\le \kappa_i^{-1}r^2)\nonumber\\ &= P_d(r^2/\kappa_i,|i|^2/\kappa_i). \label{Eq_iNotInA}
\end{align} 
Plugging \eqref{Eq_iInA} and \eqref{Eq_iNotInA} into the expression for $K(r)$ proves the claim. 
\qed }\end{proo}

\begin{rema}(The set $\mathcal{A}$ and hyperuniformity) The set $\mathcal A = \{i\in\mathbb Z^d\setminus\{0\}:\rho_i=1\}$ corresponds to perfectly correlated perturbations. Indeed, for $i\in\mathcal A$, \[ \V(p_i-p_0) = 2(1-\rho_i)\Sigma = 0, \] and therefore  $p_i-p_0=0$ \text{a.s.}

Note that under the assumptions of Corollary~\ref{Coro2}, the set $\mathcal A$ is necessarily empty. Indeed, if $\rho_i=1$ for some $i\neq0$, then $p_i-p_0=0$ almost surely and by stationarity, $p_{(k+1)i}-p_{ki}=0$ a.s. for every $k\in\mathbb N$, and therefore $p_{ki}-p_0=0$ a.s. for all $k\in\mathbb N$. Consequently, \[ C_{ki} = \C(p_{ki},p_0) = \Sigma, \] which contradicts the summability condition \[ \sum_{i\in\mathbb Z^d}\|C_i\|<\infty. \] 

\end{rema}

\begin{rema}[The case $\sigma^2=0$]
The degenerate case $\sigma^2=0$ corresponds to the stationarized lattice \[ \mathbb Z^d+U, \] for which \[ K(r) = \#\bigl((\mathbb Z^d\setminus\{0\})\cap B(0,r)\bigr), \qquad r\ge0. \] \end{rema}

\begin{rema}[Decay of the summands in $K(r)$] The terms in the series defining $K(r)$ decay exponentially fast as $||i|-r|$ increases. Indeed, for fixed $\kappa_i>0$, the quantity \[ P_d(r^2/\kappa_i,|i|^2/\kappa_i) \] is the lower tail probability of a non-central chi-squared distribution with non-centrality parameter $\kappa_i^{-1}|i|^2$, which becomes exponentially small when $|i|$ is sufficiently far from $r$. For practical simulation purposes, one is typically interested only in $r\in[0,b]$ for some fixed (and usually moderate) value $b>0$. It is therefore natural to approximate the $K$-function by the truncated sum \[ \widehat K(r) = \sum_{\substack{i\in\mathbb Z^d\setminus(\mathcal A\cup\{0\})\\ |i|\in[r-q,r+q]}} P_d(r^2/\kappa_i,|i|^2/\kappa_i) + \#(\mathcal A\cap B(0,r)), \qquad r\in[0,b], \] for a suitable truncation parameter $q>0$ depending on the dimension and the covariance structure of the perturbation field. \end{rema}


\begin{figure*}[h!]
    \centering
    \begin{subfigure}[t]{0.45\textwidth}
        \centering
        \includegraphics[height=1.2in]{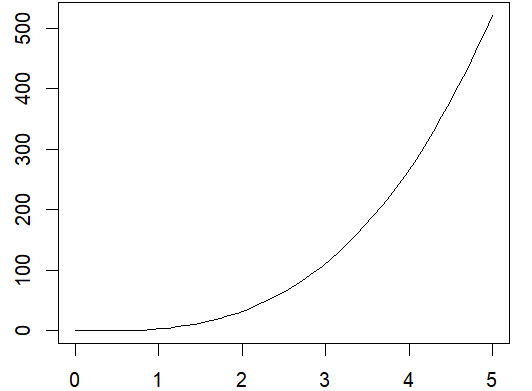}
        \caption{theoretical $K$-function}
    \end{subfigure}%
    ~ 
    \begin{subfigure}[t]{0.45\textwidth}
        \centering
        \includegraphics[height=1.2in]{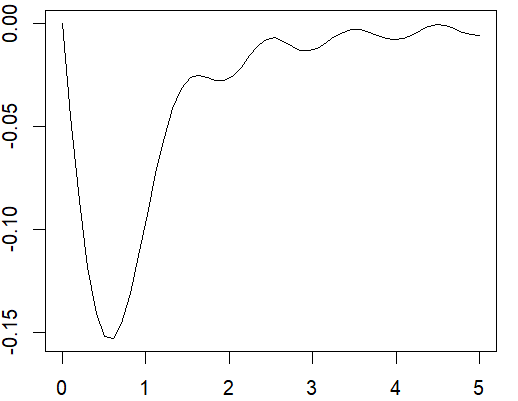}
        \caption{theoretical centered $L$-function}
    \end{subfigure}
    \caption{Theoretical values of the $K$-function and the centered $L$-function of the lattice perturbed by independent Gaussian random variables with standard deviation $0.25$.}
\end{figure*}

\section{Real data analysis}\label{S_4}
\subsection{Description of the data}

Nickel-titanium, also known as nitinol (or NiTi), is a metal alloy of nickel and titanium. Nitinol alloys are widely used in medicine, battery development, space engineering, or other fields that benefit from their recoverable strain, which can be repeatedly induced by thermal and/or mechanical actuation.

In this paper, the grain microstructure in a polycrystalline NiTi wire was reconstructed from an experimental 3D synchrotron X-ray diffraction (3D-XRD) data set. The 3D-XRD technology provides information on microstructure geometry, especially the grain center position, grain volumes, and grain orientations. However, it does not provide the full geometry of the grains, such as their exact boundaries. Reconstructing the full grain geometry, including exact grain boundaries, is the subject of ongoing numerical studies based on well-fitted mathematical models. For our study, we have a specimen at our disposal consisting of positions of 8063 grain centers which are contained in a 3D cylindrical observation window of height 80 micrometers and radius $40$ micrometers. This data set was collected by the authors in \cite{PSW19} using the so-called cross-entropy method applied to 3D X-ray diffraction measurements. To eliminate border effects in our study, we extract a 3D cubic segment of size 35x35x70 micrometers that includes 4807 grain centers.

\begin{figure}[h!]
    \centering
    \includegraphics[width=0.35\linewidth]{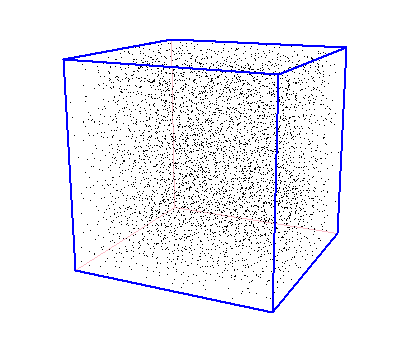}
    \caption{Positions of grain centers in NiTi alloy extracted using 3D-XRD.}
    \label{fig:enter-label}
\end{figure}

From the perspective of model fitting, this data set was extensively studied in \cite{SMB22, PSW19}. We stick with the premise of homogeneity of the data that was tested in the aforementioned papers. \textcolor{black}{The plots of the estimated summary statistics (see Figure \ref{Fig_summary} below) indicate that nearby grain centers occur less frequently than under complete spatial randomness}. 

\begin{figure}[h!]
    \centering
    \begin{subfigure}[t]{0.3\textwidth}
        \centering
        \includegraphics[height=1.2in]{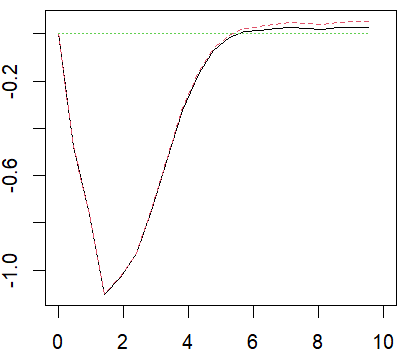}
        \caption{$\hat{L}(r)-r$}
    \end{subfigure}%
    ~ 
    \begin{subfigure}[t]{0.3\textwidth}
        \centering
        \includegraphics[height=1.2in]{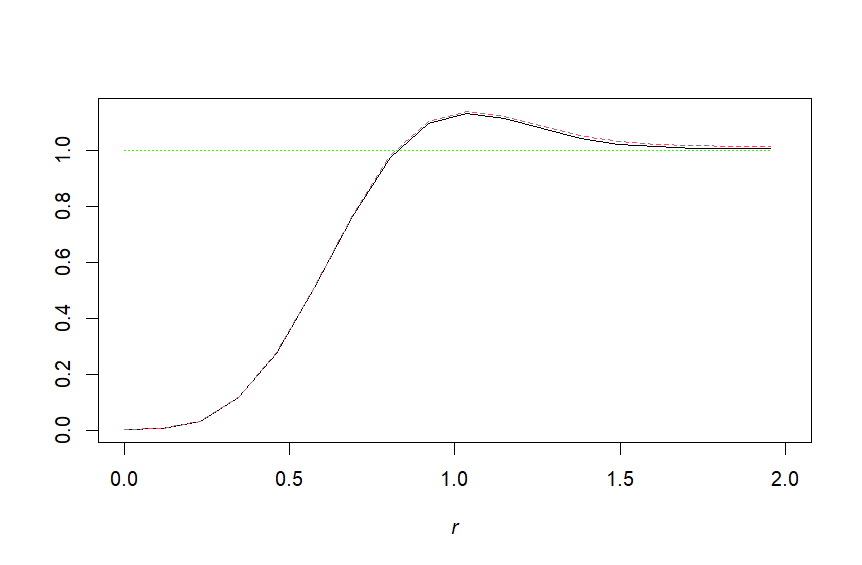}
        \caption{$\hat{g}$}
    \end{subfigure}
       ~ 
    \begin{subfigure}[t]{0.3\textwidth}
        \centering
        \includegraphics[height=1.2in]{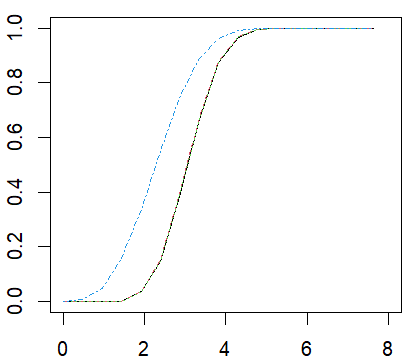}
        \caption{$\hat{G}$}
    \end{subfigure} 
    \caption{Estimated (a) centered L-function, (b) pair-correlation function and (c) G-function. The estimated values are plotted in solid lines, resp. red dashed line, while being compared to theoretical values of homogeneous Poisson point process (green dashed line).}
    \label{Fig_summary}
\end{figure}

\newpage
\subsection{Relation to the existing literature}
Primarily, we compare our approach with the methodology proposed in \cite{SMB22}. To model the positions of the grain centers in our data, the authors used the \textit{multi-Strauss process}. The finite volume model is given by a density with respect to unit rate Poisson point process
\begin{equation}\label{hierarchicalGPP}
    p_k(\mathbf{y}) = \frac{1}{Z} \beta^{\#(\mathbf{y})}\prod_{i=1}^{k-1}\gamma_i^{\sum_{y_1, y_2 \in \mathbf{y}}^{\neq} \I\{|y_1-y_2|\in (\delta_{i-1}, \delta_i]\} }, \quad \mathbf{y} \in \mathbf{N}_f,
\end{equation}
where $k\in \N$ is an a priori given parameter, $Z$ is the normalizing constant and $\beta>0, 0\leq \gamma_1, \ldots, \gamma_{k-1}\leq 1$ and $0=\delta_0<\delta_1<\ldots <\delta_{k-1}$ are unknown parameters of the model. Let $\mu_n$ be a finite volume Gibbs point process given by \eqref{hierarchicalGPP} in a window $W_n\subset \R^n$. The next lemma shows that whenever $\mu_n$ admits an accumulation point $\mu$, this is a non-hyperuniform point process in $\R^d$.

\begin{lemm}
    Any infinite volume Gibbs point process given by finite dimensional densities \eqref{hierarchicalGPP} is not hyperuniform provided that $\delta_{k-1}< \infty$.
\end{lemm}

\begin{proo}
In the setting of Section 4.1 in \cite{DF24}, the density \eqref{hierarchicalGPP} corresponds to the pairwise energy function
$$H(\mathbf{y}) = - \sum_{y_1, y_2\in \mathbf{x}}^{\neq}\sum_{i=1}^{k-1}{\I\{|y_1-y_2|\in(\delta_{i-1}, \delta_i]\}}\log \gamma_i.$$
Due to the condition posed on the parameters $\gamma_1,\ldots, \gamma_{k-1}$, the energy $H$ is non-negative regardless of the parameter $k$ and therefore, satisfies both superstability and lower-regularity assumptions. Moreover, $\mathbf{y} \in \mathbf{N}_f$ and $x \in \R^d$, the Papangelou conditional intensity is given by
$$\lambda_k^*(x, \mathbf{y})= \frac{p_k(\mathbf{y}\cup \{x\})}{p_k(\mathbf{y})}= \beta \prod_{i=1}^{k-1} \gamma_i^{\sum_{y \in \mathbf{y}}\I\{|y-x|\in (\delta_{i-1}, \delta_i]\}}.$$
Define a function $g_k:\R^d\to \R$ by
$$g_k(x)= \frac{\lambda_k^*(0, \mathbf{y}\cup \{x\})}{\lambda_k^*(0, \mathbf{y})}= \prod_{i=1}^{k-1}\gamma_i^{\I\{|x|\in(\delta_{i-1}, \delta_i]\}}.$$
The function $g_k$ does not depend on the choice of $\mathbf{y}\in \mathbf{N}_f$. With the convention $0^0=1$,
$$\int_{\R^d}\left|1- g_k(x)\right|\D x= \int_{\cup_{i=1}^{k-1}(\delta_{i-1}, \delta_i]} \left|1- g_k(x)\right|\D x\leq \sum_{i=1}^{k-1}\delta_i-\delta_{i-1}=\delta_{k-1}<\infty.$$
The resulting process is therefore not hyperuniform by Corollary 3 of \cite{DF24}.

\qed
\end{proo}

Note that this statement is in accordance with the fact that the choice $k=1$ corresponds directly to the Poisson point process with intensity $\beta$ and $k=2$ defines the Strauss process. The authors suggest modeling the data with $k=3$. That requires estimating $5$ unknown parameters of the model: $\beta, \delta_1, \delta_2, \gamma_1, \gamma_2$. In the presented methodology, we restrict ourselves to a maximum of $3$ unknown parameters. This allows us to perform faster computations compared to the latter approach.

\subsection{Hyperuniformity in the NiTi data}
\textbf{Normalizing the intensity. }The data set was first rescaled to unit intensity. The total number of points 4807 distributed in a box of volume $70^3 \mu m^3$ were resized into a box $W:= [-8.4384, 8.4384]^3$ to gain unit intensity and the center of the box being at the origin. That is in order to get us in the context of the paper \cite{MBL24} that will be used to estimate the hyperuniformity exponent and also to have the same intensity as our proposed model, a perturbed lattice. However, the latter is not necessary because the perturbed lattices can easily be defined to have any desired intensity.

\bigskip
\noindent\textbf{Isotropy}

\begin{figure}[h!]
    \centering
    \begin{subfigure}[t]{0.25\textwidth}
        \centering
        \includegraphics[height=1in]{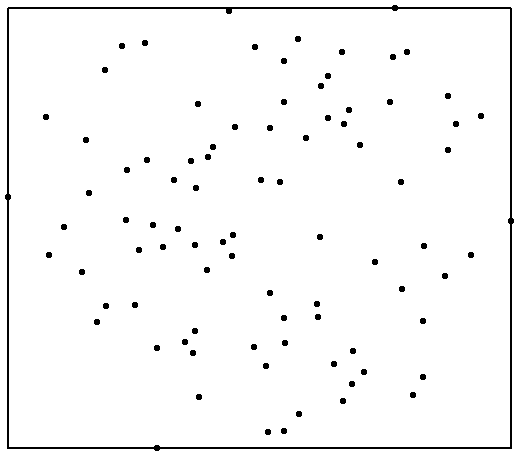}
    \end{subfigure}
       ~ 
    \begin{subfigure}[t]{0.25\textwidth}
        \centering
        \includegraphics[height=1in]{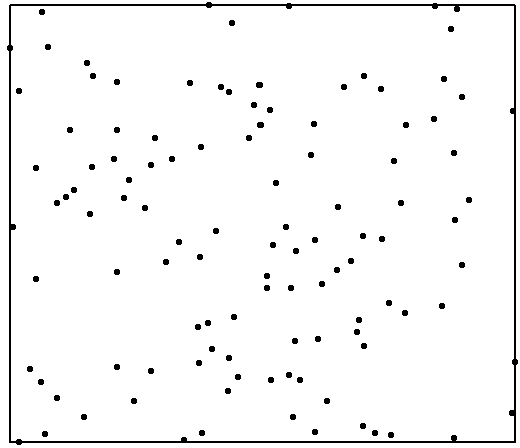}
    \end{subfigure} 
       ~ 
    \begin{subfigure}[t]{0.25\textwidth}
        \centering
        \includegraphics[height=1in]{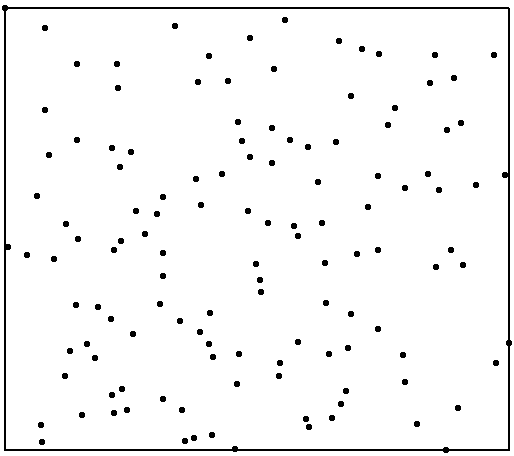}
    \end{subfigure}\\
    \begin{subfigure}[t]{0.25\textwidth}
        \centering
\includegraphics[height=0.98in]{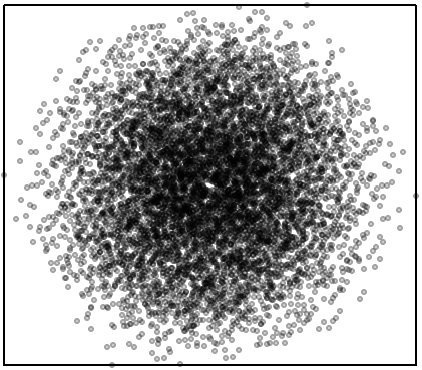}
    \end{subfigure}
        ~ 
    \begin{subfigure}[t]{0.25\textwidth}
        \centering
        \includegraphics[height=1in]{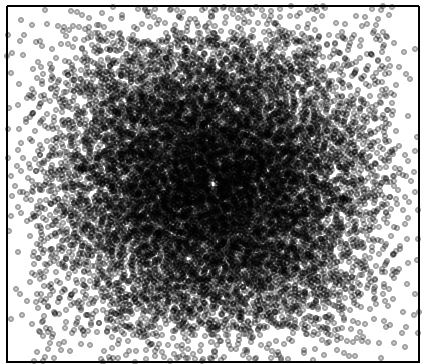}
    \end{subfigure}
        ~ 
    \begin{subfigure}[t]{0.25\textwidth}
        \centering
        \includegraphics[height=1in]{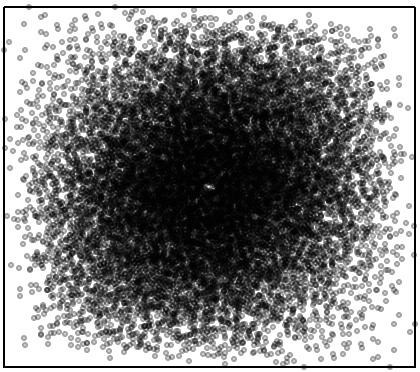}
    \end{subfigure}
    \caption{The Fry plots of the $xy$-axis slab (left), $xz$-axis slab (middle) and $yz$-axis slab (right).}\label{Fig_Fry}
\end{figure}

Possible directional effects in the data were not investigated in \cite{SMB22}.
To assess second–order isotropy of the three-dimensional point pattern, we used Fry plots (see \cite{F79}). For a point configuration $\mathbf{X}=\{x_1,\ldots, x_n\}\subset \R^3$, the Fry set is defined as the collection of all interpoint displacement vectors
$$F=\{x_i-x_j; \,x_i \neq x_j\}$$
In practice, to facilitate visualization, we consider thin planar slabs orthogonal to each coordinate axis and display the projected Fry set within each slab. Under the hypothesis that the point process is stationary and isotropic, the distribution of displacement vectors depends only on their norm $r=|x_i-x_j|$ and not on their direction. Consequently, the empirical density of points in the Fry plot is radially symmetric around the origin. Deviations from radial symmetry indicate directional dependence, which typically leads to elliptical shapes.
Fry plots are particularly informative for regular or inhibited point patterns, where preferred spacings generate visible structures in the displacement distribution \cite{RRSS18}. \textcolor{black}{This approach (see Figure \ref{Fig_Fry}) does not reveal any pronounced anisotropy at the available resolution. While some mild deviations from radial symmetry may be visually present, their magnitude appears small relative to the overall structure of the pattern. For the purposes of the subsequent modelling, we therefore proceed under the working assumption of isotropy.}

\begin{figure}
    \centering
\includegraphics[width=1\linewidth]{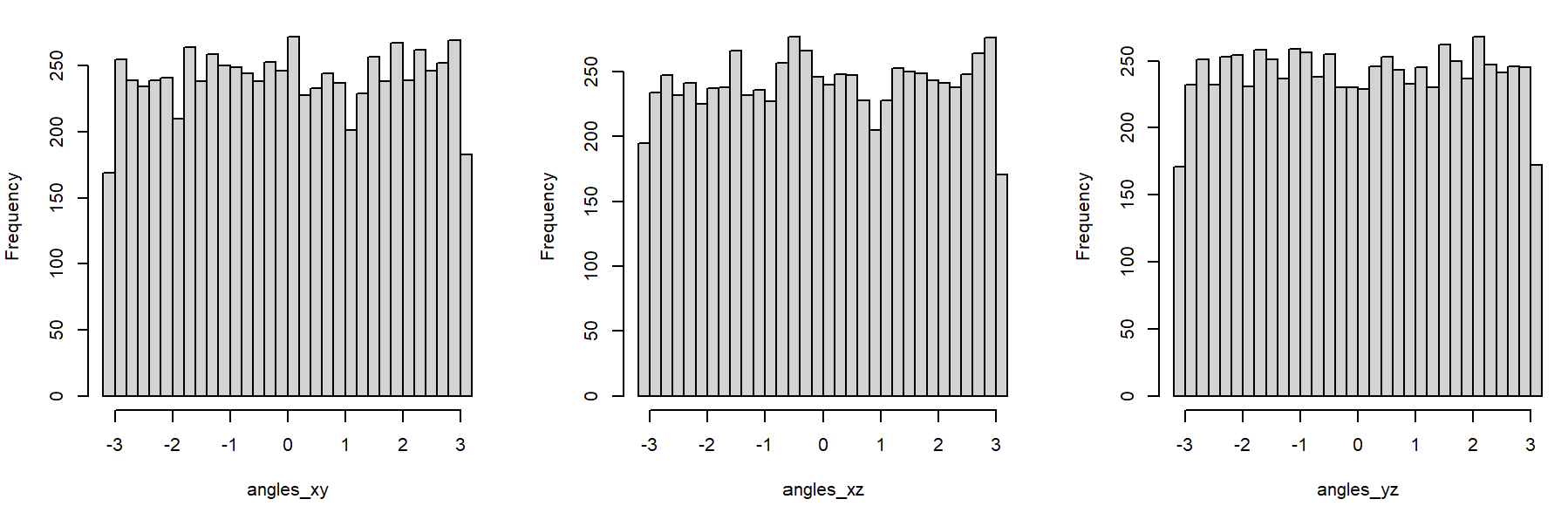}
    \caption{Empirical distribution of the angles between the points projected onto $xy$-axis (left), $xz$-axis (middle) and $yz$-axis (right) and their nearest neighbors.}
    \label{F_angles}
\end{figure}

Alternatively, we may project the points on all the axes and test the distribution of the angle between each point and its nearest neighbor. Under the hypothesis of isotropy, the distribution is uniform on $[0,2\pi]$. However, it should be noted that the nearest neighbor angles of an observed pattern are not independent. See Figure \ref{F_angles} for the empirical distribution of the angles.

\begin{figure}[h!]
    \centering
    \includegraphics[width=0.5\linewidth]{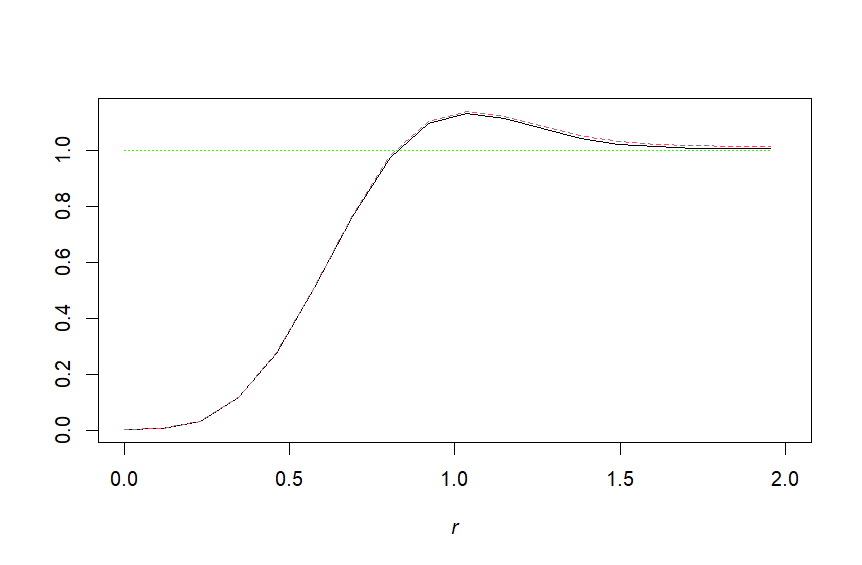}
    \caption{Estimated pair correlation function of the rescaled data set. After $r=1.5$, it becomes negligible.}
    \label{Fig_pcf}
\end{figure}

\begin{figure*}[h!]
    \centering
    \begin{subfigure}[t]{0.45\textwidth}
        \centering
        \includegraphics[height=1.5in]{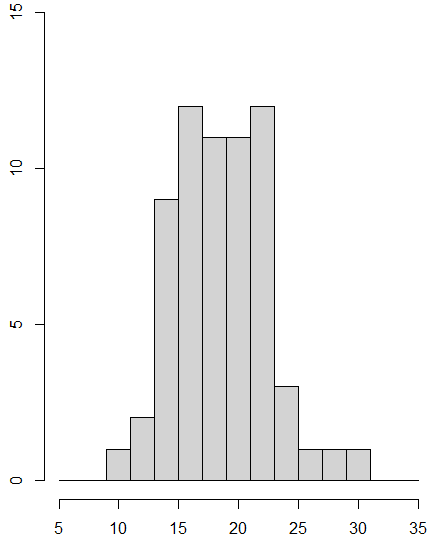}
        \caption{NiTi data}
    \end{subfigure}
       ~ 
    \begin{subfigure}[t]{0.45\textwidth}
        \centering
        \includegraphics[height=1.5in]{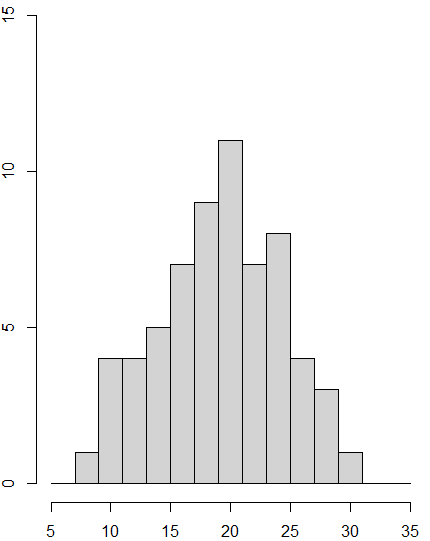}
        \caption{Poisson point process}
    \end{subfigure} 
    \caption{Histograms presenting the number of points appearing in 64 distant boxes of the size of $2.7^3$. }
    \label{Fig3}
\end{figure*}

\newpage
\bigskip
\noindent\textbf{Hyperuniformity.}
As a preliminary exploratory step, one may compare empirical fluctuations of point counts to those of a Poisson model. In Figure~\ref{Fig3}, we extracted 64 disjoint boxes separated by distance $1.5$, which corresponds to the range beyond which the estimated pair correlation function of the rescaled data becomes negligible (Figure~\ref{Fig_pcf}). The histogram for the observed data shows slightly lighter tails than that of a Poisson point process, suggesting reduced large-scale fluctuations. While this observation is consistent with hyperuniform behaviour, it is only qualitative and cannot by itself establish hyperuniformity. \textcolor{black}{In particular, the analysis is based on a single window size and therefore does not provide information about the asymptotic scaling of fluctuations, which is essential for distinguishing hyperuniformity from other fluctuation regimes.}

Statistical inference for hyperuniformity has only recently begun to develop, see in particular \cite{HGB23, KLH24, MBL24}. Some approaches require multiple independent realizations, which are not available in our setting. Detecting hyperuniformity reliably from a single finite sample remains a challenging problem, as the defining properties concern asymptotic behaviour at large scales. In dimension $d=2$, recommended sample sizes are on the order of several thousand points, and the required size is expected to increase with dimension.

Following \cite{MBL24}, we applied a test for the hyperuniformity exponent (see Eq.~\eqref{EQ_SFF}) to our dataset, obtaining the estimate $\hat{\tau}=0.81$. We provide detailed analysis in the R script at \url{https://github.com/DanielaFlimmel/Fitting_PL}. Given the limited sample size and the sensitivity of the procedure to tuning parameters and taper choice, this estimate should be interpreted cautiously as an exploratory indicator rather than a definitive asymptotic exponent. Overall, the available evidence suggests a high degree of regularity and possibly suppressed fluctuations, but it does not allow a conclusive determination of hyperuniformity. In what follows, we therefore treat hyperuniformity as a modelling hypothesis and investigate whether hyperuniform perturbed-lattice models provide an adequate description of the data.

\subsection{Models description}
We consider two models of Gaussian perturbed lattices. For these models, Proposition~\ref{L_2} provides an explicit expression for the $K$-function and can therefore be combined with the minimum contrast method for parameter estimation.

Throughout this section, $U$ denotes a random variable uniformly distributed on $[0,1]^3$ and independent of all other random variables.

\bigskip
\noindent\textbf{Independent Gaussian perturbations.} In physics papers, this model is also called an \textit{Einstein pattern}.
Assume a one-parametric model of the Gaussian perturbed lattice
\textcolor{black}{$$\Xi_{iid}^\sigma = \sum_{i \in \Z^3}\delta_{i+p^\sigma_i + U},$$}
where $\mathbf{p}=\{p^\sigma_i, i \in \Z^3\}$ are independent, identically distributed centered Gaussian random variables in $\R^3$ with variance matrix $\Sigma=\sigma^2 \mathbf{I}_3$. Clearly, $\Xi_{iid}$ forms a class-I hyperuniform point process.

\bigskip
\noindent\textbf{Dependent Gaussian perturbations.}
We consider a three-parameter family of perturbed lattices with dependent Gaussian perturbations. Specifically, we set
\textcolor{black}{
$$\Xi_{exp}^{\sigma,\delta,\gamma}= \sum_{i \in \Z^3}\delta_{i+p_i^{\sigma,\delta,\gamma}+U},$$}
where $\mathbf{p}=\{p_i^{\sigma, \delta,\gamma}, i \in \Z^3\}$ is a Gaussian random field with exponentially decaying covariances defined as follows: 
  \textcolor{black}{$$\V(p^{\sigma,\delta,\gamma}_0)= \sigma^2\mathbf{I}_3 \qquad \C(p^{\sigma,\delta,\gamma}_0, p^{\sigma,\delta,\gamma}_i) = \exp\left(
-\frac{|i|^\gamma}{\delta}
\right) \sigma^2 \mathbf{I}_3, i \in \Z^3,$$}
where $\sigma>0, \delta>0$ and $\gamma \in (0,2]$ are subjects to model fitting. The restriction $\gamma\in(0,2]$ guarantees that the covariance kernel
\begin{equation}\label{DPL_Gauss}
C(i,j):= 
\C(p_i^{\sigma,\delta,\gamma},p_j^{\sigma,\delta,\gamma})
=
\sigma^2
\exp\!\left(-\frac{|i-j|^\gamma}{\delta}\right)
I_3    
\end{equation}
 is positive definite on $\Z^3$, i.e. defines a valid stationary Gaussian random field.
Since
$$\sum_{i \in \Z^3} \|\C(p_0^{\sigma,\delta,\gamma}, p_i^{\sigma,\delta,\gamma})\|=  \sigma^2\|\mathbf{I}_3\|
\sum_{i\in\Z^3}
e^{-|i|^\gamma/\delta}<\infty,$$
the assumptions of Theorem \ref{Coro2} are satisfied. Consequently, the point process $\Xi_{exp}^{\sigma,\delta,\gamma}$ is class-I hyperuniform.


\subsection{Model fitting}
The implementation of the models together with simulation procedures is provided in the R code which is available at \url{https://github.com/DanielaFlimmel/Fitting_PL}. To eliminate the border effects of the simulations, we simulated the perturbed lattices in extended windows and then cropped it to match the size of the data.

\medskip
\noindent\textbf{Minimum contrast method. } 
The minimum contrast method is a statistical technique used to fit parametric models to observed spatial point patterns. For a rigorous introduction, see \cite{MW04}. It is particularly useful when analytic expressions for the likelihood are unavailable, but the model allows computation of theoretical summary statistics. We consider \textit{summary statistics} based on interpoint distances, i.e., functions $T:\mathcal{N}\times \R^+\to \R$, which encode structural features of the point pattern at scale $r$. Typical choices include Ripley's $K$-function, Besag's $L$-function (see \eqref{D_K} and \eqref{D_L}), or the pair correlation function, possibly after suitable transformations.

Let $\mathbf{X}$ be a realization of a point process with unknown parameter $\theta \in \Theta \subset \R^d, d\geq 1$. We denote by $T_\theta(r)$ the theoretical summary statistic under the model with parameter $\theta$, and by $\widehat T(r)$ its empirical estimator computed from the observed data $\mathbf{X}$. The minimum contrast method estimates $\theta$ by minimizing a discrepancy between $T_\theta$ and $\widehat T$ over a range of distances. 
More precisely, for a transformation $h:\mathbb{R}\to\mathbb{R}$, the contrast function is defined by
\textcolor{black}{\[
D(\theta)
=
\int_{r_{\min}}^{r_{\max}}
\bigl|h(\widehat T(r)) - h(T_\theta(r))\bigr|^2 \, dr,
\]}
where $[r_{\min},r_{\max}]$ is a prescribed interval of distances over which the summary statistic is informative and not dominated by edge effects or noise. 
The transformation $h$ (for instance $h(x)=x^{1/4}$ or $h(x)=\sqrt{x}$) is typically chosen to stabilize variance or to weight particular spatial scales.
The minimum contrast estimator is then defined as
\[
\widehat\theta
=
\arg\min_{\theta\in\Theta} D(\theta).
\]

In the case of Gaussian perturbations, we have Lemma \ref{L_2} in hand, and we can use the minimum contrast method for estimating the parameters of the models. This means to minimize
\begin{equation}\label{EQ_mincontr}
D(\theta)=\int_{r_1}^{r_2}\left|\left(\sum_{i \in \Z^3\setminus\{0\}} P_3(r^2/\kappa^\theta_i, |i|^2/\kappa^\theta_i)\right)^{1/4}-\hat{K}_{NiTi}(r)^{1/4}\right|^2\D r,
\end{equation}
where $\theta$ belongs to some parametric space given by the model and $\kappa_i^\theta$ depends on $\theta$. The estimated K-function of the NiTi data is $\hat{K}_{NiTi}$. For computational purposes, we used $\hat{K}(r)$ in \eqref{EQ_mincontr} instead of the theoretical $K$-function with $q=15$ (see the remark after Proposition \ref{lemma2}).

\medskip
\noindent\textbf{Independent Gaussian perturbations.}
For independent Gaussian perturbations, we have $\kappa_i^{\sigma}=2\sigma^2$. We choose $r_1=0, r_2=3$ and the final minimizer is $\hat{\sigma}=0.18$. Below is the comparison of the centered Besag's L-function of the fitted model together with the estimated centered L-function from the data.

\begin{figure*}[h!]
    \centering
    \begin{subfigure}[t]{0.45\textwidth}
        \centering
        \includegraphics[height=1.5in]{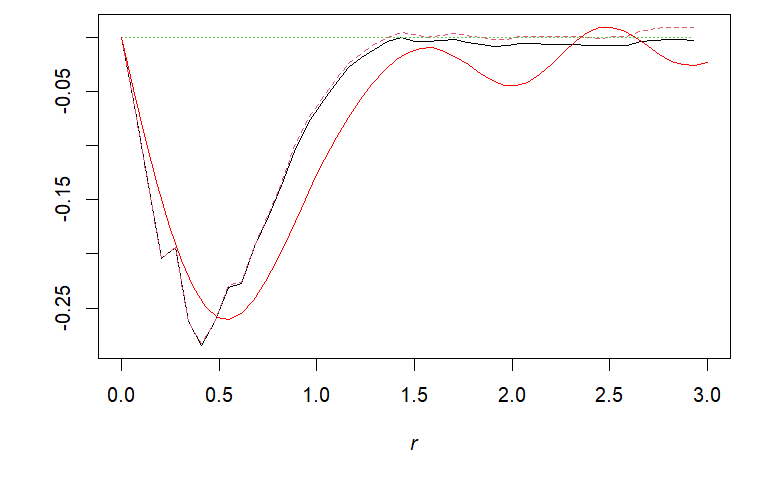}
    \end{subfigure}
       ~ 
    \begin{subfigure}[t]{0.45\textwidth}
        \centering
        \includegraphics[height=1.5in]{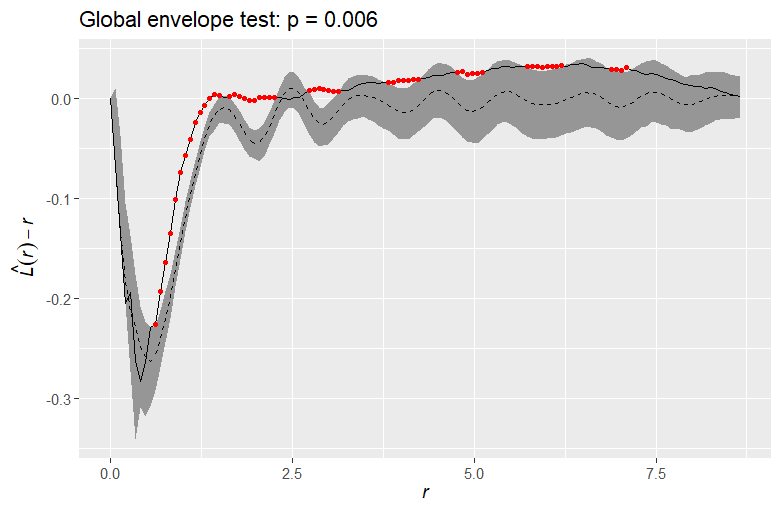}
    \end{subfigure} 
    \caption{(a) Black lines corresponds to the estimated centered L-function, green dashed line is the theoretical value of Poisson point process while red line represents the theoretical value of the fitted independently perturbed lattice with Gaussian random variables with standard deviation $\sigma=0.18$. (b) Global envelope test for the fitted i.i.d perturbed Gaussian lattice.}
    \label{Fig5}
\end{figure*}

The periodic peaks that occur at the tail of the theoretical centered $L$-function in Figure \ref{Fig5} are typical for independently perturbed lattices. Note that this phenomenon is less visible if we increase the variance in the model. However, the estimated standard deviation is low enough to preserve this property. To verify the validity of the model, we performed a global envelope test based on \cite{MMGSH17} (see Figure \ref{Fig5}, (b)). Due to these tail peaks and also suppression of the variance caused by hyperuniformity of the model, the independently perturbed lattice by Gaussian random variables is not a good approximation of the data. The need of a more flexible model motivates the subsequent paragraph.

\medskip
\noindent\textbf{Dependent Gaussian perturbations.}
Similarly as with the independent Gaussian perturbation, we used the minimum contrast method based on the $K$-function of Lemma \ref{L_2} and \eqref{EQ_mincontr} to estimate the parameters in the dependent Gaussian perturbation model resulting in $\hat{\sigma}=0.3,  \hat{\delta}=2.5$ and $\hat{\gamma} = 2$. The tail peaks were less visible in comparison to the i.i.d case. However, the global envelope test rejected this model. The reason is the low empirical variance of the $L$-function in the region $r \in (0.5, 1.5)$. This is already visible at Figure \ref{Fig5}. 

\begin{figure*}[h!]
    \centering
    \begin{subfigure}[t]{0.45\textwidth}
        \centering
        \includegraphics[height=1.5in]{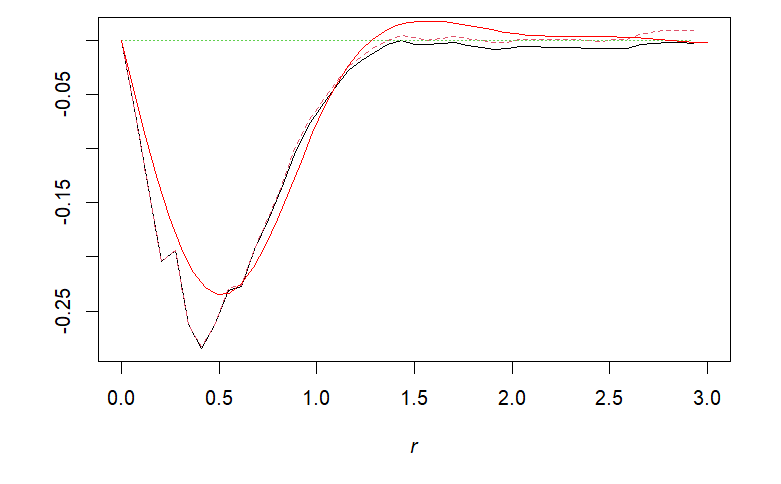}
    \end{subfigure}
       ~ 
    \begin{subfigure}[t]{0.45\textwidth}
        \centering
        \includegraphics[height=1.5in]{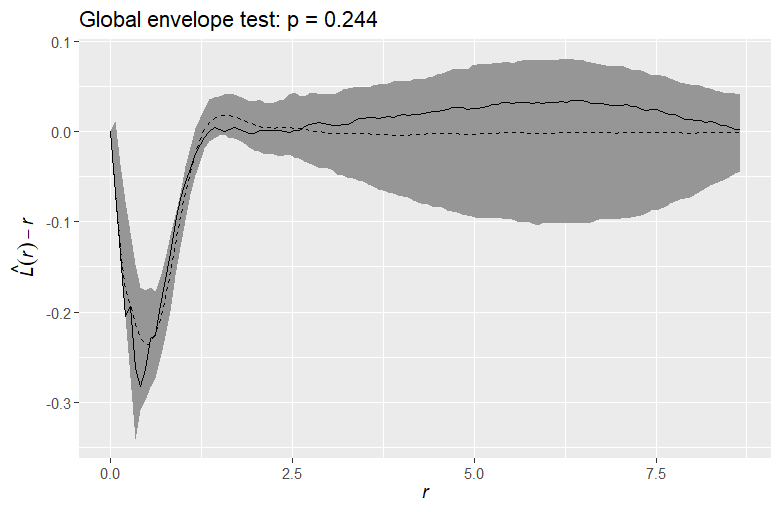}
    \end{subfigure} 
    \caption{(a) Black lines corresponds to the estimated centered L-function, green dashed line is the theoretical value of Poisson point process while red line represents the theoretical value of the fitted dependently perturbed lattice with parameters $\hat{\sigma}=0.32, \hat{\delta}=2.5$ and $\hat{\gamma}=2$ (b) Global envelope test for the fitted dependent perturbed Gaussian lattice. The corresponding p-value obtained by the global area rank envelope test is $24,4\%$}
    \label{Fig7}
\end{figure*}

For this reason, we performed the second step of the optimization procedure. Based on previous simulations, we computed the empirical variance of the centered $L$-function for each $r \in [0,3]$. Let us denote this function by $\hat{l}(r)$. The second optimization step is then based on a weighted function
\begin{equation*}
D(\theta)=\int_{r_1}^{r_2}\frac{1}{\sqrt{\hat{l}(r)}}\left|\left(\sum_{i \in \Z^3\setminus\{0\}} P_3(r^2/\kappa^\theta_i, |i|^2/\kappa^\theta_i)\right)^{1/4}-\hat{K}_{NiTi}(r)^{1/4}\right|^2\D r,
\end{equation*}
With $r_1=0.2$ and $r_2=2$, we arrive at the estimates $\hat{\sigma}=0.32,  \hat{\delta}=2.5$ and $\hat{\gamma} = 2$. The global envelope test (see Figure \ref{Fig7}) for the fitted dependent perturbed Gaussian lattice was performed resulting in $p$-value $24,4\%$ which implies a good fit.

\begin{figure*}[h!]
    \centering
    \begin{subfigure}[t]{0.45\textwidth}
        \centering
        \includegraphics[height=1.5in]{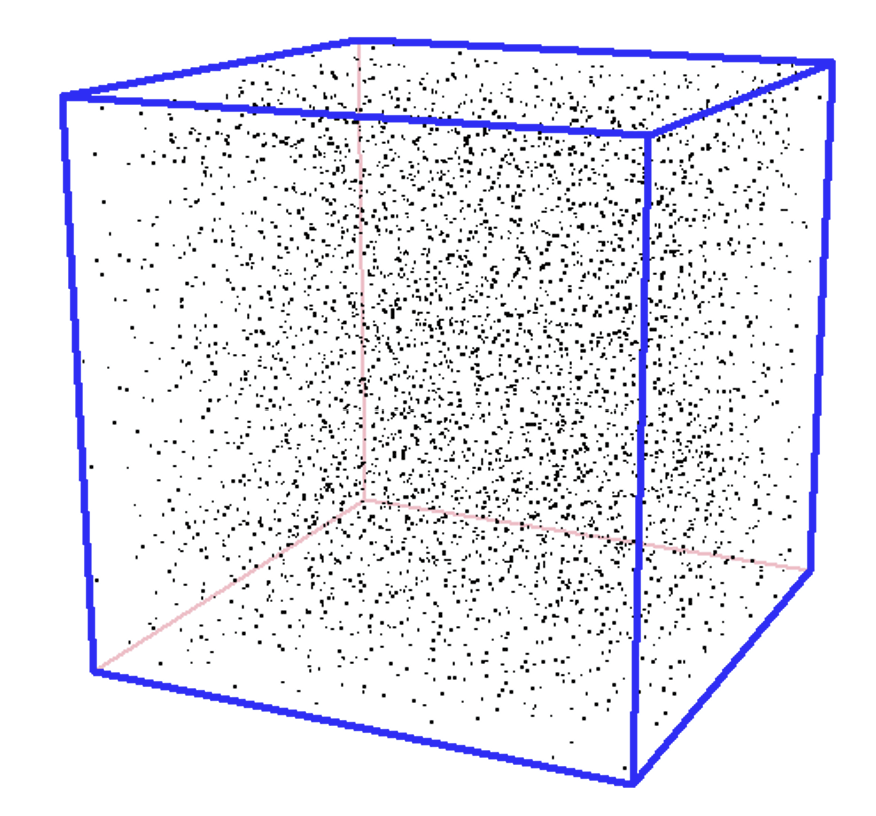}
    \end{subfigure}
       ~ 
    \begin{subfigure}[t]{0.45\textwidth}
        \centering
        \includegraphics[height=1.5in]{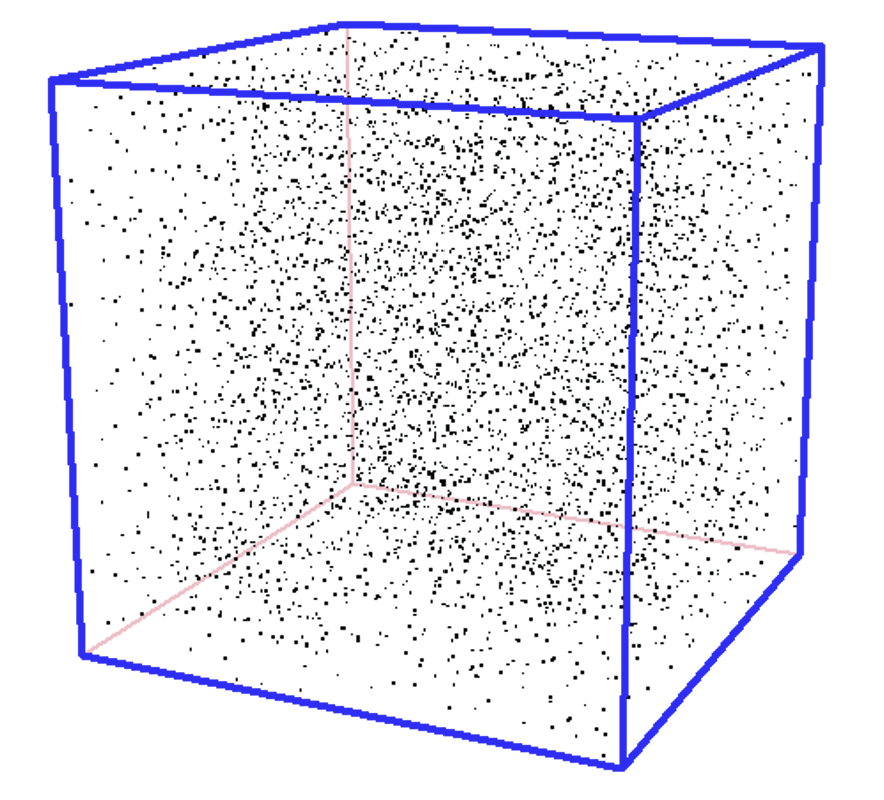}
    \end{subfigure} 
    \caption{Left is the original dataset, right is the simulated lattice perturbed by a Gaussian random field given by \eqref{DPL_Gauss} with fitted parameters $\hat{\sigma}=0.32, \hat{\delta}=2.5$ and $\hat{\gamma}=2$.}
    \label{Fig10}
\end{figure*}

Finally, we compare the original data with the simulated point pattern from the fitted model in Figure \ref{Fig10}.



\subsection{Conclusion}
When analyzing spatial data that exhibit a high degree of regularity, it is natural to investigate the growth of number variance relative to a Poisson point process. Hyperuniform models are particularly relevant when there are theoretical or physical reasons to expect suppressed large-scale fluctuations. However, reliable empirical confirmation of hyperuniformity from a single finite dataset is often challenging. In this context, perturbed lattices provide a flexible parametric framework for modelling regular point patterns, regardless of whether hyperuniformity can be shown.

From a computational perspective, models based on perturbed lattices are easy to simulate and calibrate. Independently perturbed lattices already offer a simple baseline model, although their rigid spectral structure, including visible Bragg peaks, may limit their descriptive power. Introducing correlations among the perturbations substantially increases modelling flexibility. In particular, dependent Gaussian perturbations are attractive in practice because their second-order characteristics, including the $K$-function, are available in explicit form.

Compared to classical Gibbs point processes, perturbed-lattice models avoid the need to specify and estimate an interaction potential and do not involve an intractable normalizing constant, leading to significant computational advantages. Gibbs models nevertheless remain a powerful and flexible tool for modelling data exhibiting short-range inhibition, especially when detailed control over interaction structure is required. The perturbed-lattice framework should therefore be viewed as a complementary approach that is particularly appealing when computational efficiency and analytical tractability are important considerations.

\textcolor{black}{\section{Proofs}\label{S_Proofs}}

\subsection{Proofs of the auxiliary lemmas}

The behavior of the structure factor near the origin is governed by the long-range dependence properties of the perturbation field. To quantify this dependence, we employ strong mixing coefficients. Such coefficients provide a tractable way to measure the decay of dependence between distant regions and yield bounds on covariances. These bounds will be crucial in establishing hyperuniformity. For a random variable $X$, we denote $\|X\|_p:=(\E|X|^p)^{1/p}$ for $p<\infty$ and $\|X\|_{\infty}$ the essential supremum of $|X|$ with respect to $\P$. If $X$ and $Y$ are two random variables measurable with respect to $\sigma$-fields $\mathcal{A}_1$ and $\mathcal{A}_2$ respectively, then 
two classical covariance inequalities hold (see Theorem 3 of \cite{D94} and its proof):
\begin{enumerate}
 \item If $\|X\|_\infty<\infty$, then for $q,s \geq 1$ such that $q^{-1}+s^{-1}=1$:
\begin{equation}\label{Eq_Cov1}
    |\C(X,Y)|\leq 6 \alpha^{1/s}(\mathcal{A}_1,\mathcal{A}_2)\|X\|_\infty\|Y\|_q.
\end{equation}
\item Specifically, if $X,Y$ are both a.s. bounded, then
\begin{equation}\label{Eq_Cov2}
|\C(X,Y)|\leq 4 \alpha(\mathcal{A}_1,\mathcal{A}_2)\|X\|_\infty\|Y\|_\infty,   
\end{equation}
\end{enumerate}
\textcolor{black}{where $\alpha(\mathcal{A}_1,\mathcal{A}_2)$ was defined in \eqref{Eq_mixingOfAlg}. The covariance inequalities above lead to bounds involving the two-point mixing coefficient $\alpha_{1,1}$ introduced in \eqref{Eq_mixingOfField}. More precisely, they allow us to estimate the deviation of the characteristic function of $p_i-p_0$ from its independent counterpart $|\varphi_{p_0}|^2$. This observation will be the key ingredient in the proofs of Theorem~\ref{thm1} and Corollary~\ref{thm2}.}

\begin{lemm}\label{lemma2}
    Let $\mathbf{X}=\{X_i, i \in \mathcal{L}\}$ be a strictly stationary random field such that $\E|X_0|^q<\infty$ for some $q>1$. Then there is a constant $c >0$ such that
    \begin{align*}
        \left|\varphi_{X_i-X_0}(t)-|\varphi_{X_0}(t)|^2\right| & \leq c \min\{ \textcolor{black}{\alpha_{1,1}^{1/s}}(|i|)|t|\left(\E|X_0|^q\right)^{1/q}, \textcolor{black}{\alpha_{1,1}}(|i|)\} \quad \text{for all } t \in \R^d,     \end{align*}
   where $1/s+1/q=1$. If, moreover, $X_0$ is centered and $\E|X_0|^{2q}<\infty$, then there is a constant $c'>0$ such that
   \begin{align*}
        \left|\varphi_{X_i-X_0}(t)-|\varphi_{X_0}(t)|^2\right| & \leq c' \min\{ \textcolor{black}{\alpha^{1/s}_{1,1}}(|i|)|t|^2\left(\E |X_0|^{2q}\right)^{1/q},\textcolor{black}{\alpha_{1,1}}(|i|)\} \quad \text{for all } t \in \R^d.     \end{align*}
   Here, $\varphi_Y$ denotes the characteristic function of the random vector $Y$.
\end{lemm}

\begin{proo}
First, assume that $t \in B(0,1)$. In a similar spirit as in the proof of Theorem 2.1. of \cite{LSC11}, we write
\begin{align*}
    \left|\varphi_{X_k-X_0}(t)-|\varphi_{X_0}(t)|^2\right| & = \left|\E e^{i t \cdot(X_k-X_0)}-\E e^{i t \cdot X_k}\E e^{-i t \cdot X_0}\right|\\
    &= \left| \E \cos(t\cdot(X_k-X_0))+i \E \sin(t(X_k-X_0)) \right.\\
    & \quad \left.- \left(\E\cos(t \cdot X_k)+ i \E \sin(t\cdot X_k)\right)\left(\E\cos(t \cdot X_0) - i \E \sin(t\cdot X_0)\right)\right|\\
\end{align*}
Using the relations $\sin(x+y) = \sin(x)\cos(y)+\cos(x)\sin(y),$ $\cos(x+y)=\cos(x)\cos(y)-\sin(x)\sin(y)$, $\sin(x)=-\sin(-x)$ and $\cos(x)=\cos(-x)$, the latter expression can be estimated by $I_1 + I_2 + I_3 + I_4$, where
\begin{align*}
I_1 & := | \C(\cos(t\cdot X_k), \cos(t\cdot X_0))|,\\
I_2 &:= |\C(\sin(t\cdot X_k), \sin(t\cdot X_0))|,\\
I_3 & := |\C(\sin(t\cdot X_k), \cos(t\cdot X_0))|,\\
I_4 & :=|\C(\cos(t\cdot X_k), \sin(t\cdot X_0))|. 
\end{align*}
For $I_2, I_3, I_4$, we note that
$$\max\{\E |\sin(t\cdot X_k)|^q,\E |\sin(t\cdot X_0)|^q \} \leq |t|^q \E|X_0|^q<\infty.$$
On the other hand, $\max\{\E|\cos(t\cdot X_k)|,\E|\cos(t\cdot X_0)|\}\leq 1$. Using the covariance inequality \eqref{Eq_Cov1} for strong mixing random fields,
$$\max \{I_2,I_3, I_4\} \leq 6 |t | \| X_0 \| _q \textcolor{black}{\alpha_{1,1}^{1/s}}(|k|).$$
At last, having $\cos(x)= 1- 2 \sin^2(x/2)$ and $\sin^2(x)\leq |\sin(x)|\leq |x|$, we arrive at
\begin{align*}
    I_1 &=|\C(\cos(t\cdot X_k), 1- 2\sin^2(t\cdot X_0/2))|\\
    & = 2|\C(\cos(t\cdot X_k),\sin^2(t\cdot X_0/2))|\\
    & \leq c |t| \|X\|_q \textcolor{black}{\alpha_{1,1}^{1/s}}(|k|). 
\end{align*}

This bound can be improved if $\E X_0=0$ and $\E|X_0|^{2q}<\infty$. We observe that $\E \sin(t X_0)=\E \sin(t X_k)=0$ and so the initial estimate reduces to  
$$\left|\varphi_{X_k-X_0}(t)-|\varphi_{X_0}(t)|^2\right| \leq I_1.$$
Then, using the bound $\sin^2(x)\leq x^2$ gives
$$I_1=2|\C(\cos(t\cdot X_k),\sin^2(t\cdot X_0/2))| \leq c' |t|^2 \left(\E |X|^{2q}\right)^{1/q} \textcolor{black}{\alpha_{1,1}^{1/s}}(|k|).$$

For $|t|>1$, we simply use the covariance inequality \eqref{Eq_Cov2} for bounded random variables $X= e^{i t \cdot X_k}, Y= e^{-i t \cdot X_0}$ to write
 $$\left|\varphi_{X_k-X_0}(t)-|\varphi_{X_0}(t)|^2\right| = |\C(Y,X)|\leq 4 \textcolor{black}{\alpha_{1,1}}(|k|).$$

\qed
\end{proo}

To assess the speed of convergence of $\langle S, j_r\rangle$, we make use of the following asymptotic results for the kernels $j_r$. 
\begin{lemm}\label{L_bessel}
    For any $d\geq 1$,
\begin{align}
    \int_{B(0,1)}|x| j_r(x) \D x & = O\left( \frac{\log r}{r}\right),\label{EQ_jrx}\\
    \int_{B(0,1)}|x|^2 j_r(x) \D x & =O\left(\frac{1}{r}\right)\label{EQ_jrx2}\\
    \int_{B(0,1)^C}j_r(x)\D x &= O\left(\frac{1}{r}\right) \label{EQ_jrBc}.
\end{align}
\end{lemm}
\begin{proo}
To omit the constants in what follows, we will write $f(r) \asymp g(r)$ if the functions have the same asymptotic order, i.e.
$$0 < \lim_{r \to \infty} \frac{f(r)}{g(r)}<\infty.$$
This is the same as writing $f(r) = O(g(r))$.
 First, we evaluate the integral \eqref{EQ_jrx}. The integrand is a radial function, so the substitutions $y=|x|$ and consequently, $z = r y$ lead to
    \begin{align*}
        \int_{B{(0,1)}}|x| j_r(x) \D x & \asymp \int_{B(0,1)} \frac{J_{d/2}(r|x|)^2}{|x|^{d-1}}\D x\\
        &\asymp \int_0^1 J_{d/2}(r y) ^2\D y\\
        & = \int_0^{r}\frac{J_{d/2}(z)^2}{r}\D z  
    \end{align*}
    Next, we use the classical Bessel asymptotics (cf \cite[Section 5.16]{LS72}):
    $$J_{d/2}(z) = O(z^{d/2}) \text{ for }z \to 0,$$
    \begin{equation}\label{EQ_BesselInf}
        J_{d/2}(z) = O\left(\frac{1}{\sqrt{z}}\right) \text{ for }z \to \infty.
    \end{equation}
 With $r >1$, we have
\begin{align*}
    \int_0^{r}\frac{J_{d/2}(z)^2}{r}\D z & = \int_0^{1}\frac{J_{d/2}(z)^2}{r}\D z +\int_1^{r}\frac{J_{d/2}(z)^2}{r}\D z\\
    & \asymp\frac{1}{r}\int_0^1 z^d \D z + \frac{1}{r}\int_1^{ r} \frac{1}{z}\D z\\
    & \asymp \frac{\log r}{r}
    .
\end{align*}
The proofs of \eqref{EQ_jrx2} and \eqref{EQ_jrBc} are analogous, leading to
$$\int_{B(0,1)}|x|^2 j_r(x) \D x \asymp \frac{1}{r^2}\int_0^1 z^{d+1} \D z + \frac{1}{r^2}\int_1^r 1 \D x = O\left(\frac{1}{r}\right),$$
$$\int_{B(0,1)^C}j_r(x)\D x \asymp\int_r^{\infty}\frac{J_{d/2}(z)^2}{z}\D z \asymp \int_r^{\infty}\frac{1}{z^2}\D z = O\left(\frac{1}{r}\right).$$
\qed
\end{proo}

\subsection{Hyperuniformity of the stationarized and independently perturbed lattices}
The result in Lemma \ref{L_statL} dates back to the \textit{Gauss circle problem}. In general, counting integer points in Euclidean balls is a difficult task. However, viewing the problem through the spectral representation and calculus with distributions leads to a short and elegant proof. 

In the following, we make use of the Poisson summation formula that can be found in the present form in \cite[Chapter VII.]{SW71}.
\begin{lemm}[Poisson summation formula]
    Let $f:\R^d \to \mathbb{C}$ be integrable and there exist $C, \delta>0$ for which $\max\{|f(x)|, |\hat{f}(x)|\}\leq \frac{C}{|x|^{d+\delta}}$ as $|x|\to \infty$. Then
    $$\sum_{x \in \mathcal{L}}f(x)= \sum_{k \in \mathcal{L^*}}\hat{f}(2\pi k),$$
    where $\mathcal{L}^*:= \{k \in \R^d \text{ such that } k \cdot x \in \Z \text{ for all } x \in \mathcal{L}\}$ is the dual lattice to $\mathcal{L}$.
\end{lemm}

\begin{proo}[of Lemma \ref{L_statL}]
    The reduced second factorial cumulant measure of $\Xi$ is
    $$\gamma_{red}^{(2)}(B)=\#(\mathcal{L}\setminus\{0\}\cap B)-\lambda^d(B)$$
 and we may rewrite \eqref{EQ_generalSF} as 
\begin{equation}\label{EQ_L1_1}
\langle S, j_r\rangle=1+\sum_{k \in \mathcal{L}\setminus\{0\}}\hat{j_r}(k)-\int \hat{j_r}(x) \D x.
\end{equation}
By definition of $j_r$, we have
$$\hat{j_r}(0) = \frac{\I_{B_r}\star \I_{B_r}(0)}{\lambda^d(B_r)}=\frac{1}{\lambda^d(B_r)}\int \I_{B_r}(x) \I_{B_r}(0-x) \D x = 1$$
and
\begin{equation}\label{EQ_jr0}
\int \hat{j_r}(x) \D x = \langle \lambda^d, \hat{j_r}\rangle = \langle \hat{\lambda^d}, j_r \rangle = \langle (2\pi)^d \delta_0, j_r\rangle = (2\pi)^d j_r(0).
\end{equation}
In the last line, we used the fact that $\hat{\lambda^d}= (2\pi)^d \delta_0$. We plug these expressions into \eqref{EQ_L1_1} and use the Poisson summation formula. Then
\begin{align*}
    \langle S, j_r\rangle & = \sum_{k \in \mathcal{L}}\hat{j_r}(k)-(2\pi)^d j_r(0)\\
     & = \sum_{x \in \mathcal{L}^*}\hat{\hat{j_r}}(2\pi x)-(2\pi)^d j_r(0)\\
    & = (2\pi)^d\sum_{x \in \mathcal{L}^*}j_r(-2\pi x)-(2\pi)^d j_r(0)=(2\pi)^d\sum_{x \in \mathcal{L}^*\setminus\{0\}}j_r(2\pi x).
\end{align*}
To prove class-I hyperuniformity, we use the asymptotic property \eqref{EQ_BesselInf}. As a consequence, we have $j_r(x) = O(r^{-1})$ for any $x \neq 0.$ Moreover, there are constants $C, r_0>0$ such that $|j_r(x)| \leq C \frac{1}{|x|^{d+1}}$ for every $x \neq 0$ and $r>r_0$. As a result, 
$$\langle S, j_r\rangle = O(r^{-1}).$$
    \qed
\end{proo} 

\begin{proo}[of Lemma \ref{L_IPL}]
\textcolor{black}{For independently perturbed lattices, the moment assumption on $p_0$ from Proposition 2.6 in \cite{DFHL24} is not needed for local square integrability. Indeed, for any bounded Borel set $B$, Tonelli's theorem gives
$$
\E \Xi(B)=\E \sum_{i\in\mathcal L}\I\{i+p_i+U\in B\} = \E \sum_{i\in\mathcal L}\I\{i+p_0+U\in B\}.
$$
Since $B$ is bounded, there exists $n\in \N$ and disjoint half-open unit cubes $Q_1, \ldots,Q_n$ with integer vertices such that $B \subset \cup Q_k$. For any $z \in \R^d$ and $k=1,\ldots ,n$ there exists a unique lattice point $i_k+z \in Q_k$. Taking $z=p_0+U$ gives
$$\E\Xi(B)\leq n.$$
By independence, 
\begin{align*}
\E \Xi(B)^2 & = \E \sum_{i \in \mathcal{L}}\I\{i + p_i + U \in B\}+ \E \sum_{\substack{i,j \in \mathcal{L}\\i\neq j}}\I\{i+p_i+U \in B\}\I\{j+p_j+U \in B\} \\
& \leq n + n^2<\infty.
\end{align*}
Thus the second-order measures are well defined without imposing any moment assumption on $p_0$ and using \eqref{alpha2PSL}, we may write
    $$\gamma_{red}^{(2)}(B) = \E \sum_{k \in \mathcal{L}\setminus\{0\}}\I\{k+p_k-p_0 \in B\}-\lambda^d(B).$$}

Next, we make use of the non-negativity of $\hat{j_r}$, \eqref{EQ_jr0} and identity in distribution of $p_k - p_0$, for all $k \in \mathcal{L}\setminus\{0\}$:

\begin{align*}
    \langle S, j_r\rangle & = 1 + \E \sum_{k \in \mathcal{L}\setminus \{0\}} \hat{j_r}(k + p_k - p_0) - \int \hat{j_r}(x) \D x \\
    & = 1+ \sum_{k \in \mathcal{L}\setminus \{0\}} \E \hat{j_r}(k + p_k - p_0) - (2\pi)^d j_r(0)\\
    & =   1-\E \hat{j_r}(p_1-p_0)+ \sum_{k \in \mathcal{L}} \E \hat{j_r}(k + p_1-p_0) - (2\pi)^d j_r(0)\\
    & = 1- \int |\varphi_{p_0} (x)|^2 j_r(x) \D x +  \E \sum_{k \in \mathcal{L}} \mathcal{F}\left[e^{-i (p_1-p_0)\cdot} j_r(\cdot)\right] (k)- (2\pi)^d j_r(0)\\
    & = 1- \int |\varphi_{p_0} (x)|^2 j_r(x) \D x +  \E \sum_{x \in \mathcal{L}^*} \mathcal{F}\left[\mathcal{F}\left[e^{-i (p_1-p_0)\cdot} j_r(\cdot)\right]\right] (2\pi x)- (2\pi)^d j_r(0)\\
    & = 1- \int |\varphi_{p_0} (x)|^2 j_r(x) \D x +  (2\pi)^d\E \sum_{x \in \mathcal{L}^*} e^{i (p_1-p_0)2\pi x} j_r(-2\pi x) - (2\pi)^d j_r(0)\\
     & = 1- \int |\varphi_{p_0} (x)|^2 j_r(x) \D x +  (2\pi)^d \sum_{x \in \mathcal{L}^*\setminus\{0\}} |\varphi_{p_0}(2\pi x)|^2 j_r(2\pi x).
\end{align*}
It remains to prove the hyperuniformity. First, the sum $\sum_{x \in \mathcal{L}^*\setminus\{0\}} |\varphi_{p_0}(2\pi x)|^2 j_r(2\pi x)$ is dominated by 
 $\sum_{x \in \mathcal{L}^*\setminus\{0\}} j_r(2\pi x)$ which is $O(r^{-1})$ by Lemma \ref{L_statL}. Next, the function $g(x):=|\varphi_{p_0}(x)|^2$ is continuous, non-negative, bounded by $1$ and $g(0)=1$. Therefore, for every $\varepsilon>0$ there exists $\delta>0$ such that
$$1-\varepsilon\leq g(x)\leq 1 + \varepsilon, \quad \text{for all } x \in B_\delta.$$ 
Consequently,
$$(1-\varepsilon)\int_{B_\delta}j_r(x) \D x \leq \int_{B_\delta} g (x) j_r(x) \D x \leq (1+\varepsilon)\int_{B_\delta}j_r(x) \D x.$$
On the other hand
$$0\leq \int_{B^C_\delta}g(x)j_r(x)\D x \leq \int_{B_\delta ^C}j_r(x)\D x.$$
Since, $j_r\to \delta_0$ and $\varepsilon$ can be arbitrarily small, we conclude that
$$\lim_{r\to \infty}\int |\varphi_{p_0} (x)|^2 j_r(x) \D x = 1$$
which shows hyperuniformity. 
\qed
\end{proo}

\subsection{Proofs of the main results}
\begin{proo}[of Theorem \ref{thm1}]
 We follow the ideas of the proof of Proposition 3.1 of \cite{DFHL24}. In fact, the assumption on $\alpha$-mixing coefficient enables drastic simplifications of the respective proof. We briefly recall the necessary notation first. The moment assumption on the perturbations guarantees that $\E[\textcolor{black}{\Xi(B_r)^2}]<\infty$ for all $r>0$ and the reduced second factorial cumulant measure is a well-defined locally finite Borel measure given by
 \begin{equation*}
     \gamma^{(2)}_{red} (B)= \E \sum_{k \in \mathcal{L}\setminus \{0\}}\I\{k+p_k-p_0 \in B\}- \lambda^d(B), \quad B \in \mathcal{B}^d_b.
 \end{equation*}
 In combination with \eqref{EQ_generalSF} we arrive at
 \begin{equation*}
 \sigma(r)= \frac{\V (\textcolor{black}{\Xi(B_r)})}{\lambda^d(B_r)}=1+ \E \sum_{k \in \mathcal{L}\setminus \{0\}}\hat{j_r}(k+p_k-p_0 )- \int \hat{j_r}(x)\D x, 
 \end{equation*}
 where $j_r$ is the smooth approximation of the mass around $0$ discussed in Section \ref{S_GeneralSF}. The strategy to show that $\sigma(r)\to 0$ is to first, extract the contribution of the stationarized lattice and second, compare the rest with independently perturbed lattice (cf. Definition \ref{D_PerturbedLattice}). For further reference, we denote the stationarized lattice by $\Xi_{stat}$ and the independently perturbed lattice with perturbations following the law of $p_0$ by $\Xi_{iid}$.
 The second step reduces to comparing characteristic functions $\varphi_{p_k-p_0}$ of the random variable $p_k-p_0$ with function $|\varphi_{p_0}|^2$ which is the characteristic function of $p_k-p_0$ when $p_k$ and $p_0$ are independent and identically distributed.

As proposed above, let us split $\sigma(r)$ so that
 $$\sigma(r)= A(r)+B(r),$$
 where
 \begin{align*}
     A(r)& = 1+\sum_{k \in \mathcal{L}\setminus \{0\}} \hat{j_r}(k)-\int\hat{j_r}(x)\D x,\\
     B(r)& = \E  \sum_{k \in \mathcal{L}\setminus \{0\}} \hat{j_r}(k+p_k-p_0) - \hat{j_r}(k).
 \end{align*} 
 
Recall that $A(r)=\V(\textcolor{black}{\Xi_{stat}(B_r)})/\lambda^d(B_r)$ and by Lemma \ref{L_statL}, it is $O(r^{-1})$.
The residual term $B(r)$ can be further expressed using the characteristic functions $\varphi_{p_k-p_0}, k \in \mathcal{L}\setminus\{0\}$:
$$
    B(r)=\sum_{k \in \mathcal{L}\setminus\{0\}}\int_{\R^d} e^{-i t \cdot k} j_r(t)(\varphi_{p_k-p_0}(t)-1)\D t.
$$
As the mass of $j_r$ concentrates around the origin, it is convenient to further denote 
\begin{align*}
    B_1(r)&:=\sum_{k \in \mathcal{L}\setminus\{0\}}\int_{B(0,1)} e^{-i t \cdot k} j_r(t)(\varphi_{p_k-p_0}(t)-1)\D t \\
    B_2(r)&:=\sum_{k \in \mathcal{L}\setminus\{0\}}\int_{B(0,1)^C} e^{-i t \cdot k}j_r(t)(\varphi_{p_k-p_0}(t)-1)\D t
\end{align*}
so that $B(r)=B_1(r)+B_2(r).$
For $\Xi_{iid}$, we denote the latter quantities by $B^{iid}(r), B^{iid}_1(r)$, resp. $B^{iid}_2(r)$. We know that $B^{iid}(r) \to 0$ as $r \to \infty$ as a consequence of Lemma  \ref{L_IPL}.  

First, using Lemma \ref{lemma2} with the assumption $\|p_0\|_q:=(\E|p_0|^q)^{1/q}<\infty$,
\begin{align*}
    |B_1(r)-B^{iid}_1(r)| & = \left| \sum_{k \in \mathcal{L}\setminus\{0\}}\int_{B(0,1)} e^{i t \cdot k}j_r(t)\left(\varphi_{p_k-p_0}(t)-|\varphi_{p_0}(t)|^2\right) \D t\right|\\
    & \leq \sum_{k \in \mathcal{L}\setminus\{0\}}\int_{B(0,1)} j_r(t) \left| \varphi_{p_k-p_0}(t)-|\varphi_{p_0}(t)|^2\right| \D t\\
    & \leq c \|p_0\|_q \int_{B(0,1)} j_r(t)|t| \D t \sum_{k \in \mathcal{L} \setminus \{0\}} \textcolor{black}{\alpha_{1,1}^{1/s}}(|k|).
\end{align*}
This term is $O(\log(r)/r)$ by \eqref{EQ_jrx}.

Similarly, again by Lemma \ref{lemma2},
\begin{align*}
    |B_2(r)-B^{iid}_2(r)| & \leq  \sum_{k \in \mathcal{L}\setminus\{0\}}\int_{B(0,1)^C} j_r(t) \left| \varphi_{p_k-p_0}(t)-|\varphi_{p_0}(t)|^2\right| \D t\\
    & \leq c \int_{B(0,1)^C} j_r(t) \D t \sum_{k \in \mathcal{L} \setminus \{0\}} \textcolor{black}{\alpha_{1,1}}(|k|).
\end{align*}
Since $\sum_{k \in \mathcal{L}}\alpha^{1/s}(|k|)$ converges, so does the sum $\sum_{k \in \mathcal{L}\setminus\{0\}}\alpha(|k|).$ Eventually, this term is $O(r^{-1})$ by \eqref{EQ_jrBc}.

If $p_0$ is symmetric (again without loss of generality around 0) and $\E|p_0|^{2q}<\infty$, then Lemma \ref{lemma2} gives an improved bound
$$|B_1(r)-B^{iid}_1(r)|\leq c' \left(\E|p_0|^{2q}\right)^{1/q}\int_{B(0,1)}j_r(t) |t|^2\D t \sum_{k \in \mathcal{L} \setminus \{0\}} \textcolor{black}{\alpha_{1,1}^{1/s}}(|k|),$$
which is $O(r^{-1})$ by \eqref{EQ_jrx2}.

\qed
\end{proo}

The proof of Corollary \ref{thm2} is analogous. In the proof of Lemma \ref{lemma2}, we would use two times the covariance inequality \eqref{Eq_Cov2} to arrive at the respective estimate 
 \begin{align*}
        \left|\varphi_{X_k-X_0}(t)-|\varphi_{X_0}(t)|^2\right| & \leq c \min\{ \textcolor{black}{\alpha_{1,1}}(|k|)|t|\|X_0\|_{\infty}, \textcolor{black}{\alpha_{1,1}}(|k|)\} \quad \text{for all } t \in \R^d.   \end{align*}

\begin{proo}[of Theorem \ref{Coro2}]
\textcolor{black}{For $d=1$, the assertion follows from \cite{DFHL24}, since Gaussian perturbations possess moments of all orders. The following proof works for general dimension. Essentially, the argument follows the proof of Theorem \ref{thm1} almost verbatim, replacing the mixing bound by an explicit Gaussian estimate.
Let 
\[
\Sigma:=\V(p_0),
\qquad
C_k:=\C(p_k,p_0).
\]
For $t\in\mathbb R^d$, Gaussianity and stationarity give
\[
\varphi_{p_k-p_0}(t)
=
\exp\left(
-\frac12 t^\top\V(p_k-p_0)t
\right)
=
\exp\left(
-\frac12 t^\top(2\Sigma-C_k-C_k^\top)t
\right),
\]
and also
\[
|\varphi_{p_0}(t)|^2=\exp(-t^\top\Sigma t).
\]
Consequently,
\[
\left|
\varphi_{p_k-p_0}(t)-|\varphi_{p_0}(t)|^2
\right|
=
e^{-t^\top\Sigma t}
\left|
\exp\left(\frac12t^\top(C_k+C_k^\top)t\right)-1
\right|.
\]
By the Cauchy--Schwarz inequality,
\[
|t^\top C_k t|
=
|\C(t^\top p_k,t^\top p_0)|
\leq
\sqrt{\V(t^\top p_k)}
\sqrt{\V(t^\top p_0)}
=
t^\top\Sigma t.
\]
Consequently, setting
\(
x_k(t):=\frac12t^\top(C_k+C_k^\top)t,
\)
\begin{equation}\label{Eq_allk}
|x_k(t)|
\leq
t^\top\Sigma t.
\end{equation}
Using the elementary inequality
\[
|e^x-1|\le |x|e^{|x|},
\]
it follows from \eqref{Eq_allk} that
\[
e^{-t^\top\Sigma t}|e^{x_k(t)}-1|
\le
|x_k(t)|e^{-t^\top\Sigma t+|x_k(t)|}
\le
|x_k(t)|.
\]
Therefore,
\begin{equation}\label{EQ_GaussianBasic}
\left|
\varphi_{p_k-p_0}(t)-|\varphi_{p_0}(t)|^2
\right|
\leq
|x_k(t)|.
\end{equation}}

\textcolor{black}{
We now consider the region $|t|\leq1$. Since
\[
|x_k(t)|
\leq
\frac12\|C_k+C_k^\top\|\,|t|^2
\leq
\|C_k\|\,|t|^2,
\]
the estimate \eqref{EQ_GaussianBasic} yields
\begin{equation}\label{EQ_GaussianEstimate}
\left|
\varphi_{p_k-p_0}(t)-|\varphi_{p_0}(t)|^2
\right|
\leq
\|C_k\|\,|t|^2.
\end{equation}
Adopting the notation from the proof of Theorem~\ref{thm1} and using \eqref{EQ_GaussianEstimate} yields
\begin{align*}
|B_1(r)-B_1^{iid}(r)| & \leq  \sum_{k\in\mathbb Z^d\setminus\{0\}}
\int_{B(0,1)}
j_r(t)
\left|
\varphi_{p_k-p_0}(t)-|\varphi_{p_0}(t)|^2
\right|
\D t\\
&\leq \sum_{k\in\mathbb Z^d\setminus\{0\}}\|C_k\|
\int_{\mathbb R^d} j_r(t)|t|^2\,\D t.
\end{align*}
By the summability assumption and \eqref{EQ_jrx2}, the latter is $O(r^{-1})$.}

\textcolor{black}{
It remains to estimate the contribution from $B(0,1)^C$. Let $\lambda$ denote the smallest eigenvalue of $\Sigma$. Since $\Sigma$ is positive definite, $\lambda>0$ and
\[
t^\top \Sigma t \geq \lambda |t|^2,
\qquad t\in\mathbb R^d.
\]
Since $\sum_{k\in\mathbb Z^d}\|C_k\|<\infty,$ there exists $k_0 \geq 0$ such that $\|C_k\|\leq \frac{1}{2}\lambda$ for all $|k|> k_0$. For such $k$ we have
\begin{equation}\label{EQ_GaussImproved}
|x_k(t)|
\leq
\|C_k\||t|^2
\leq
\frac12\lambda|t|^2
\leq
\frac12 t^\top\Sigma t.
\end{equation}
Using once again the inequality $|e^x-1|\leq |x|e^{|x|}$ and \eqref{EQ_GaussImproved}, we obtain
\begin{align*}
\left|
\varphi_{p_k-p_0}(t)-|\varphi_{p_0}(t)|^2
\right|
&=
e^{-t^\top\Sigma t}|e^{x_k(t)}-1|\\
&\leq
|x_k(t)|e^{-t^\top\Sigma t+|x_k(t)|}\\
&\leq
\|C_k\||t|^2
e^{-t^\top\Sigma t/2}\\
&\leq
\|C_k\||t|^2
e^{-\lambda|t|^2/2}.
\end{align*}
Therefore,
\begin{align*}
&\sum_{\substack{k \in \mathcal{L}\\ |k|> k_0}}
\int_{B(0,1)^C}
j_r(t)
\left|
\varphi_{p_k-p_0}(t)-|\varphi_{p_0}(t)|^2
\right|
\D t \leq 
\sum_{\substack{k \in \mathcal{L}\\ |k|> k_0}}
\|C_k\|
\int_{B(0,1)^C}
j_r(t)|t|^2
e^{-\lambda|t|^2/2}\D t.
\end{align*}
Since
\(
|t|^2e^{-\lambda|t|^2/2}
\)
is bounded on $B(0,1)^C$, the integral above is bounded by a constant multiple of
\(
\int_{B(0,1)^C}j_r(t)\,dt,
\)
which is $O(r^{-1})$ by \eqref{EQ_jrBc}. For $|k| \leq k_0$, we use the trivial bound
\(
\left|
\varphi_{p_k-p_0}(t)-|\varphi_{p_0}(t)|^2
\right|
\leq 2,
\) and \eqref{EQ_jrBc} to obtain for some constant $C_0<\infty$
\begin{align*}
&\sum_{\substack{k\in \mathcal{L}\setminus
\{0\}\\ |k|\leq k_0}}
\int_{B(0,1)^C}
j_r(t)
\left|
\varphi_{p_k-p_0}(t)-|\varphi_{p_0}(t)|^2
\right|
\D t
\leq  C_0
\int_{B(0,1)^C}j_r(t)\,\D t
=
O(r^{-1}),
\end{align*}
Eventually,
$$|B_1(r)-B_1^{iid}(r)|+|B_2(r)-B_2^{iid}(r)| = O(r^{-1}),$$
which finishes the proof.}
\qed
\end{proo}

\section{Future extension}


\noindent\textbf{Central limit theorems for cell characteristics of a tessellation based on the NiTi data.}
First, we should find a suitable model for the distribution of the radii $r$ (e.g. by following the hierarchical approach by \cite{SMB22}) to construct a marked perturbed lattice $\Xi_m=\{(x,m), x \in \Xi, r\geq 0\}=\{(x_i,r_i)_{i \in \Z^3}\}$, where the cell generated by a marked point $(x,r)$ is
$$C^\rho((x,r), \Xi_m)= \bigcap_{i \in \Z^3}\{y \in \R^3: \rho(y,(x,r))\leq \rho(y,(x_i,r_i))\}.$$
The common choices of the function $\rho$ are listed below
\begin{itemize}
    \item Voronoi tessellation: $\rho(y,(x,r))= |y-x|$
    \item Laguerre tessellation: $\rho(y,(x,r))= |y-x|^2-r^2$
    \item Johnson-Mehl tessellation: $\rho(y,(x,r))= |y-x|-r$
\end{itemize}
If $K_\mathbf{0}^\rho(\Xi_m)$ is the cell generated
by the typical point of $\Xi_m$ (in the Palm sense), we are interested in the estimation of $\E \xi(K_\mathbf{0}^\rho)$ where $\xi$ is some function on closed sets (diameter, volume, etc.)
$$H_n= \sum_{i \in \Z^3\cap W_n} \xi(C((x_i,r_i),\Xi_m)) \cdot \text{weight}_i,$$
where the weights are chosen reasonably so that the estimator $H_n$ is (asymptotically) unbiased and consistent. The stabilization methods and mixing properties of the perturbation field could then be used to show that
$$\frac{H_n-\E H_n}{\sqrt{\sigma_n}}\xrightarrow[n\to \infty]{\mathcal{D}} N(0,\sigma^2_\xi),$$
where $\sigma_n= \textbf{var}|\Xi\cap W_n|$.
This asymptotic study should be carried out in a separate paper.

\section*{Acknowledgement}
This work was supported by the Czech Science Foundation, project no. 22-15763S.

I am grateful to Gabriel Mastrilli for his assistance with the implementation of the test for the hyperuniformity exponent in dimension $d=3$, to Jiří Dvořák for his valuable suggestions concerning the model-fitting part, and to Luca Lotz for insightful observations that led to improvements of the section devoted to Gaussian perturbations.

I also thank the anonymous referees for their careful reading and constructive comments, which resulted in substantial improvements to both the content and presentation of the paper.

Finally, I am deeply grateful to my husband and parents-in-law for their unwavering support and for helping with childcare responsibilities, which made this research possible


\begin{thebibliography}{99}
\footnotesize

\bibitem{BH24} Björklund, M. and Hartnick, T.  (2024): \textit{Hyperuniformity and non-hyperuniformity of quasicrystals}, Math. Ann. 389, 365--426.





\bibitem{C21}
Coste, S. (2021): \textit{Order, fluctuations, rigidities.} Available \hyperlink{https://scoste.fr/assets/survey_hyperuniformity.pdf}{here}.


\bibitem{DVJ08}
Daley, D. J. and Vere-Jones, D. (2003): \textit{An Introduction to the
Theory of Point Processes: Volume I: Elementary Theory and Methods},
Second Edition. Springer, New York.


\bibitem{DFHL24}
Dereudre, D., Flimmel, D. Huesmann, M. and Leblé, T.  (2024): \textit{(Non)-hyperuniformity of perturbed lattices.} Preprint. \hyperlink{https://arxiv.org/pdf/2405.19881}{arxiv.org/pdf/2405.19881}. 

\bibitem{DF24}
Dereudre, D. and Flimmel, D. (2024): \textit{Non-hyperuniformity of Gibbs point processes with short-range interactions}, J. Appl. Probab. \textbf{61} (4), 1380--1406. 

\bibitem{D94}
Doukhan, P. (1994): \textit{Mixing: Properties and Examples}, Lecture Notes in Statistics vol. 85, Springer-Verlag. 

\bibitem{F79} Fry, N. (1979): \textit{Random point distributions and strain measurement in rocks.} Tectonophysics, vol. 60, 89--105.

\bibitem{G04}
Gabrielli, A. (2004): \textit{Point processes and stochastic displacement fields}, Phys. Rev. E. \textbf{70}(6), 066131.

\bibitem{GS75}
Gács, P and Szász, D. (1975): \textit{On a Problem of Cox Concerning Point Processes in $\R^k$ of "Controlled Variability"}. Ann. Probab. \textbf{3}(4), 597--607.


\bibitem{GL17}
Ghosh, S. and Lebowitz, J. (2017):
\textit{Fluctuations, large deviations and rigidity in hyperuniform systems: a brief survey}, Indian J. Pure Appl. Math. \textbf{48}(4), 609--631.


\bibitem{HGB23}
Hawat, D., Gautier, G., Bardenet, R. et al. (2023): \textit{On estimating the structure factor of a point process, with applications to hyperuniformity}, Stat. Comput. \textbf{33} (61). 

\bibitem{H95} 
Hof, A. (1995): \emph{Diffraction by aperiodic structures at high temperatures}, J. Phys. A 28, 57--62.

\bibitem{JS25}
Jalowy, J and Stange, H. (2025): \textit{Box-covariances of hyperuniform point processes.} Preprint. \hyperlink{https://arxiv.org/abs/2506.13661}{arxiv.org/abs/2506.13661}.



\bibitem{JLH14}
Jiao, Y., Lau, T. and Hatzikirou, H. et al (2014): \textit{Avian photoreceptor patterns represent a disordered hyperuniform solution to a multiscale packing problem}, Phys. Rev. E 89. 

\bibitem{K75}
Kanter, M. (1975): \textit{Stable Densities Under Change of Scale and Total Variation Inequalities}, Ann. Probab. \textbf{3}(4), 697--707.


\bibitem{KKT20}
Klatt, M. A., Kim, J. and Torquato, S. (2020): \textit{Cloaking the underlying long-range order of randomly perturbed lattices}, Phys. Rev. E 101.



\bibitem{KLH24}
Klatt, M. A., Last, G. and Henze, N. (2024): \textit{A genuine test for hyperuniformity.} Preprint. \hyperlink{https://arxiv.org/abs/2210.12790}{arxiv.org/abs/2210.12790}.

\bibitem{KLLY25}
Klatt, M. A., Last, G., Lotz, L. and Yogeshwaran D. (2025): \textit{Invariant transports of stationary random measures: asymptotic variance, hyperuniformity, and examples.
} Preprint. \hyperlink{https://arxiv.org/abs/2506.05907}{arxiv.org/pdf/2506.05907}

\bibitem{LRY24}
Lachièze-Rey, R. and Yogeshwaran, D. (2024): \textit{Hyperuniformity and optimal transport of point processes.} Preprint. \hyperlink{https://arxiv.org/abs/2402.13705}{arxiv.org/abs/2402.13705}.


\bibitem{LR25}
Lachièze-Rey, R. (2025): \textit{Hyperuniform random measures, transport and rigidity}. Preprint. \hyperlink{https://arxiv.org/abs/2510.18392}{arxiv.org/abs/2510.18392}.

\bibitem{LP17}
Last, G. and Penrose, M. (2017): \textit{Lectures on the Poisson process.} Cambridge University Press, Cambridge.

\bibitem{LS72}
Lebedev, N.N. and Silverman, R.A. (1972): \textit{Special Functions and Their Applications.} Dover Books on Mathematics. Dover Publications.


\bibitem{LSC11}
Li, Y., Shanchao, Y. and Chengdong, W. (2011): \textit{Some inequalities for strong mixing random variables with applications
to density estimation}, Stat. Probab. Lett. 81, 250--258.

\bibitem{MBL24}
Mastrilli, G., B\l{}aszczyszyn, B. and Lavancier, F. (2024): \textit{Estimating the hyperuniformity exponent of point processes.} Preprint. \hyperlink{https://arxiv.org/pdf/2407.16797}{arxiv.org/pdf/2407.16797}.

\bibitem{MBMW15}
Mayer, A., Balasubramanian, V., Mora, T. and Walczak, A. M. (2015): \textit{How a well-adapted immune system is organized}, Proc. Nat. Acad. Sci. E 112, 5950--5955.

\bibitem{MW04}
Møller, J. and Waagepetersen, R. P. (2004): \textit{Statistical Inference and Simulation for Spatial Point Processes}, 2nd edition, Chapman \& Hall/CRC.

\bibitem{MMGSH17}
Myllymäki, M., Mrkvička, T.,  Grabarnik, P., Seijo, H. and Hahn, U. (2017): \textit{Global envelope tests for spatial processes.} J. R. Stat \textbf{79}(2), 381--404


\bibitem{P93}
Peebles, P. J. E. (1993): \textit{Principles of Physical Cosmology}, Princeton University Press, Princeton. 


\bibitem{PS14}
Peres, Y. and Sly, A. (2014): \textit{Rigidity and tolerance for perturbed lattices.} Preprint. \hyperlink{https://arxiv.org/abs/1409.4490}{arxiv.org/pdf/1409.4490}. 

\bibitem{PSW19}
Petrich, L., Staněk, J., Wang, M. et al. (2019): \textit{Reconstruction of Grains in Polycrystalline Materials From Incomplete Data Using Laguerre Tessellations}, Microsc. Microanal. \textbf{25}(3), 743--752. 


\bibitem{RRSS18}
Rajala, T., Redenbach, C., Särkkä, A. and Sormani, M. (2018): \textit{A review on anisotropy analysis of spatial point patterns.} Spatial Statistics 28, 141--168.

\bibitem{SPY07}
Schreiber, T., Penrose, M. D. and Yukich, J. E. (2007): 
\textit{Gaussian Limits for Multidimensional Random Sequential Packing at Saturation.} 
Commun. Math. Phys. 272, 167--183.


\bibitem{SMB22}
Seitl, F., Møller, J. and Beneš, V. (2022): \textit{Fitting three-dimensional Laguerre tessellations by hierarchical marked point process models}, Spat. Stat. 51. 

\bibitem{ST06}
Sodin, M. and Tsirelson, B. \textit{Random complex zeroes, II. Perturbed lattice.} Isr. J. Math. 152, 105--124.

\bibitem{SW71}
Stein, E. M. and Weiss, G. (1971): \textit{Introduction to Fourier analysis on Euclidean spaces}, Volume 1. Princeton university press.

\bibitem{TS03}
Torquato, S. and Stillinger, F. H. (2003): \textit{Local density fluctuations, hyperuniformity, and order metrics}, Phys. Rev. E 68, 041113.

\bibitem{T16}
Torquato, S. (2016): \textit{Hyperuniformity and its generalizations}, Phys. Rev. E 94, 022122.

\bibitem{T18}
Torquato, S. (2018): \textit{Hyperuniform states of matter}, Phys. Rep. 745, 1--95. 

\bibitem{Y21} Yakir, O. (2021): \textit{Fluctuations of linear statistics for Gaussian perturbations of the lattice $\mathbb{Z}^d$}, J. Stat. Phys. 182, Paper No.~58.


\bibitem{Y22} Yakir, O. (2022): \textit{Recovering the lattice from its random perturbations}, Int. Math. Res. Not.
IMRN 8, 6243--6261.
\end{thebibliography}
\end{document}